\documentclass[11pt]{amsart}

\newcommand{\private}[1]{}
\usepackage{amssymb, amsmath, amscd, amsthm, color, epsfig, url, graphics}

\usepackage[all]{xy}          
\xyoption{dvips}              

\makeatletter
\renewcommand\l@subsection{\@tocline{2}{0pt}{2pc}{5pc}{}}
\makeatother

\newcommand{\R}{{\mathbb R}}

\newcommand{\Q}{{\mathbb Q}}
\newcommand{\abs}[1]{{\left\vert #1 \right\vert}}

\newcommand{\hofiber}{\operatorname{hofiber}}

\newcommand{\holim}{\operatorname{holim}}
\newcommand{\colim}{\operatorname{colim}}
\newcommand{\hocolim}{\operatorname{hocolim}}
\newcommand{\tfiber}{\operatorname{tfiber}}

\newcommand{\Map}{\operatorname{Map}}

\newcommand{\Emb}{\operatorname{Emb}}

\newcommand{\Imm}{\operatorname{Imm}}
\newcommand{\Link}{\operatorname{Link}}

\newcommand{\interior}{\operatorname{int}}

\newcommand{\del}{{\partial}}

\newcommand{\Top}{\operatorname{Top}}

\theoremstyle{plain}
\newtheorem{thm}{Theorem}[section]
\newtheorem{prop}[thm]{Proposition}

\newtheorem{cor}[thm]{Corollary}

\theoremstyle{definition}
\newtheorem{defin}[thm]{Definition}
\newtheorem{example}[thm]{Example}
\newtheorem{def/ex}[thm]{Definition/Example}

\theoremstyle{remark}
\newtheorem{rem}[thm]{Remark}

\newcommand{\refS}[1]{Section~\ref{S:#1}}
\newcommand{\refT}[1]{Theorem~\ref{T:#1}}
\newcommand{\refC}[1]{Corollary~\ref{C:#1}}
\newcommand{\refP}[1]{Proposition~\ref{P:#1}}
\newcommand{\refD}[1]{Definition~\ref{D:#1}}

\newcommand{\refE}[1]{equation~$(\ref{E:#1})$}

\begin{document}


\title[Introduction to the manifold calculus of Goodwillie-Weiss]{Introduction to the manifold calculus of Goodwillie-Weiss}


\author{Brian A. Munson}
\address{Department of Mathematics, Wellesley College, Wellesley, MA}
\email{bmunson@wellesley.edu}
\urladdr{http://palmer.wellesley.edu/\~{}munson}





\begin{abstract}
We present an introduction to the manifold calculus of functors, due to Goodwillie and Weiss. Our perspective focuses on the role the derivatives of a functor $F$ play in this theory, and the analogies with ordinary calculus. We survey the construction of polynomial functors, the classification of homogeneous functors, and results regarding convergence of the Taylor tower. We sprinkle examples throughout, and pay special attention to spaces of smooth embeddings.
\end{abstract}

\maketitle

\tableofcontents

\parskip=4pt
\parindent=0cm

\section{Introduction}\label{S:intro}

We intend to explain some of the intuition behind one incarnation of calculus of functors, namely the so-called ``manifold calculus'' due to Weiss and Goodwillie \cite{W:EI1,GW:EI2}. Specifically, we will highlight some analogies between the ordinary calculus of functions $f:\R\to\R$ and the manifold calculus of functors. The trouble with analogies is that they are not equivalences, and some may lead the reader to want to push them further. Some may indeed be pushed further than we are currently aware, and some may lead to direct contradictions and/or bad intuition. Another risk is that it is considered bad manners to tell people how to categorize various ideas: part of our mathematical culture seems to be that we leave intuition for talks and personal communications and rigor and precision for our papers, and with good reason: we cannot anticipate the ways in which our work may be useful in the future, and so it may be best to convey it in as concise and precise a way as possible. We feel the relatively small risk of misleading the reader and the faux pas of making permanent intuitive notions by publishing them is a small price to pay for the possibility that this may entice some reader to learn more about these ideas and try to use them. Finally, we would like to emphasize that this is not meant to be a rigorous introduction to calculus of functors. We will frequently omit arguments which we find distract us from our attempts to be lighthearted. We hope this work makes digesting the actual details from the original sources easier for newcomers.

The philosophy of calculus of functors is to take a functor $F$ and replace it by its Taylor series, and we will begin our discussion of ordinary calculus there and work backwards. Associated to a smooth function $f:\R\to\R$ is its Taylor series at zero (we choose zero for convenience; any center will work just fine):
\begin{equation}\label{E:ordtaylorseries} 
f(0)+f'(0)x+f''(0)\frac{x^2}{2!}+\cdots+f^{(n)}(0)\frac{x^n}{n!}+\cdots.
\end{equation}
There are two natural questions to ask about this power series: (1) does it converge, and if so, for what $x$?; and (2) if it converges, does it converge to $f$? The Taylor series is computationally much easier to work with than the function. A typical application is to truncate the series at degree $k$, thus obtaining the \emph{$k^{th}$ degree Taylor polynomial} $T_kf$ of the function $f$. If one is lucky and $f^{(k+1)}(x)$ can be controlled to be small in some neighborhood of zero, then one can use Taylor's inequality to estimate the remainder. Specifically, if $|f^{(k+1)}(x)|\leq M$ in a neighborhood of zero, then the remainder
$$
R_k(x)=|f(x)-T_k(x)|\leq M\frac{|x|^{k+1}}{(k+1)!}.
$$
Our first goal is to make the analog of a Taylor series for a functor $F$ which associates each open set $U$ in a smooth manifold $M$ a topological space (we will specify the categories and hypotheses on $F$ soon). A simple example to keep in mind is the space of maps, $U\mapsto\Map(U,X)$ for some space $X$. Our second goal will be to explore issues of convergence in the special case of spaces of embeddings; here the functor of interest is $U\mapsto\Emb(U,N)$, the space of smooth embeddings of $U$ in a smooth manifold $N$. Here are some natural questions that arise based on the above:
\begin{enumerate}
\item What is the definition of the derivative of a functor, and how should we compute it?
\item What is the definition of a polynomial functor?
\item Can we write the Taylor polynomial of a functor as polynomial whose coefficients are the derivatives?
\item What is a ``good'' approximation, and what should we mean by convergence?
\end{enumerate}

We will answer these questions and more. Despite our attempt at lightheartedness, there will be no avoiding certain constructions unpleasant for the purposes of an introductory paper. The main culprit here is homotopy limits and colimits, and we will assume the reader is more or less familiar with these. If the reader has not seen these before or has only a nodding acquaintance with them, let her not despair; we will try to give some intuition about what role these objects play, though it may still remain largely indigestible. Nevertheless, we have done two things: (1) Provided an appendix with the statements and attempts at explanation of results we have used in proofs, and (2) We have tried give alternate, hopefully simpler, constructions whenever possible, and focused on special cases where intimate knowledge of homotopy limits and colimits is not necessary. We also assume the reader is familiar with the basics of differential topology, namely handlebody decompositions of smooth manifolds, and the basics of transversality.

This paper is organized into two main parts. The first, \refS{prelims} to \refS{homogclass}, is concerned with developing the notion of derivatives and polynomials, and tell us how to build a Taylor series for a functor, and \refT{homogclass} even gives a reasonable description of its homogeneous pieces. The usefulness of the definitions developed in these sections pay off in the proof of \refT{charlinear}, and this proof contains a useful organizational principle important to later arguments, namely induction on the handle dimension. The second part is devoted to the question of convergence. We make a few general remarks about convergence in \refS{conv}, and then move on to the specific case of spaces of embeddings in \refS{convemb}.

\subsection{Further reading}\label{S:more}

This work is an introduction, not a sample of the state-of-the art, but it is right of the reader to ask whether there is any point to this endeavor, so what follows are some references to applications of manifold calculus. Any omissions are due to the ignorance of the author. Weiss himself wrote a more rigorous survey \cite{W:EmbCalc} with a different perspective than this one. A survey with emphasis on spaces of embeddings and related spaces is \cite{GKW:survey}. Another survey with many ideas from differential topology which are useful for studying spaces of embeddings using manifold calculus is \cite{Bud:FES}. For a survey on homotopy calculus, which shares many of the same tools as manifold calculus, see \cite{Kuhn:GT&CH}. As for applications of manifold calculus to spaces of embeddings, there are several recent works, including  \cite{ALTV}, \cite{Arone:DerI}, \cite{ALV}, \cite{Bud:LCLK}, \cite{BCSS}, \cite{GK:SE}, \cite{GK:PE}, \cite{GKW}, \cite{LTV:Vass}, \cite{M:Emb}, \cite{MV:Links}, \cite{MV:Multi}, \cite{ScaSinha}, \cite{S:TSK}, \cite{S:OKS},  \cite{V:IT}, and \cite{V:FTK}. For applications to spaces of link maps and connections with generalizations of Milnor's invariants, see \cite{GM:LinkStable}, \cite{M:Milnor}, and \cite{M:LinkNumber}.

\subsection{Conventions}\label{S:convs}

We will not be too careful about the category of spaces in which we will work. For some purposes, the category $\Top$ of compactly generated spaces will be fine. For other purposes, such as spaces of maps, we work in the category of simplicial spaces (a $k$-simplex in $\Map(X,Y)$ is a map $\Delta^k\to\Map(X,Y)$). We will, by abuse, always let $\Top$ be the target category. We write $\underline{k}$ in place of $\{1,2,\ldots, k\}$. We let $\interior(X)$ stand for the interior of a subset $X$ of some topological space. A space $X$ is $k$-connected if $\pi_i(X)$ vanishes for $0\leq i\leq k$ for all choices of basepoint in $X$. Every space is $(-2)$-connected, and nonempty spaces are $(-1)$-connected. A map of spaces $f:X\to Y$ is $k$-connected if it is an isomorphism on $\pi_i$ for $1\leq i <k$ and a surjection when $i=k$ for all possible choice of basepoints. Its homotopy fibers are therefore $(k-1)$-connected. Conversely, if for all choice of basepoints in $Y$, the homotopy fiber of $f$ is $(k-1)$-conected, then $f$ is $k$-connected. In particular, every map is $(-1)$-connected.

The union of a smooth manifold $L^m$ with boundary with a ``$j$-handle'' $H^j=D^j\times D^{n-j}$ is obtained by choosing an embedding $e:S^{j-1}\to \del L$ and forming the identification space $L\cup_f H^j$ by attaching $H^j$ to $\del L$ along $\del D^j\times D^{n-j}\subset H^j$. We refer to $j$ as the \emph{dimension} of the handle, and refer to $D^j\times\{0\}$ as the \emph{core} of the handle. All smooth compact manifolds admit a ``handle decomposition'', which is a description of $M$ in terms of handles of various dimensions and attaching maps (see \cite{Kos:DiffMflds}). We define the \emph{handle dimension of $M$} to be the smallest integer $j$ such that $M$ admits a handle decomposition with handles of dimension less than or equal to $j$. A handlebody decomposition of a smooth manifold $M^m$ is the analog of a cell structure on $M$, with $j$-handle playing the role of $j$-cell. Note that if $P^p$ is a smooth compact submanifold of $M$, then the disk bundle of its normal bundle is a smooth compact codimension $0$ submanifold of $M$ of handle dimension at most $p$.

A $j$-dimensional handle $H^j$ is a manifold with corners. That is, $\del H^j=\del D^j\times D^{n-j}\cup D^j\times\del D^{n-j}$, and this union happens along the corner set $\del D^j\times\del D^{n-j}$. More generally we will eventually encounter what is called a \emph{smooth manifold triad}. Roughly speaking, this is a triple $(Q,\del_0Q, \del_1Q)$, which is a smooth manifold $Q$ of dimension $q$ whose boundary is decomposed as $\del_0Q\cup\del_1Q$ and whose corner set is $\del_0Q\cap\del_1Q$.  Boundary points have neighborhoods which look locally like $[0, \infty)\times \R^{q-1}$, and points in the corner set have neighborhoods which look locally like $[0,\infty)\times[0,\infty)\times \R^{q-2}$. In particular, we regard a $j$-handle $H^j$ as a smooth manifold triad with $\del_0H^j=\del D^j\times D^{n-j}$ and $\del_1H^j=D^j\times\del D^{n-j}$. We refer the reader to \cite{GW:EI2} for details.

\section{Preliminaries}\label{S:prelims}

We need to discuss the axioms necessary to impose on our functors to obtain an interesting and computable theory.

\begin{defin}
Let $M$ be a smooth closed manifold of dimension $m$. Define $\mathcal{O}(M)$ to be the category (poset) of open subsets of $M$. Its objects are open sets $U\subset M$, and morphisms $U\to V$ are the inclusion maps $U\subset V$.
\end{defin}

Manifold calculus studies contravariant functors $F:\mathcal{O}(M)\to\Top$ which satisfy two axioms. Before we state them, let us consider a few examples, all of which are basically some sort of space of maps.

\begin{example}
Let $X$ be any space. The functor $\Map(-,X):\mathcal{O}(M)\to\Top$ given by the assignment $U\mapsto\Map(U,X)$ is a contravariant functor, since an inclusion $U\subset V$ gives rise to a restriction map $\Map(V,X)\to\Map(U,X)$.
\end{example}

\begin{example}
Let $N$ be a smooth manifold. The embedding functor $\Emb(-,N):\mathcal{O}(M)\to\Top$ is given by $U\mapsto\Emb(U,N)$. This is the space of smooth maps $f:U\to N$ such that (1) $f$ is one-to-one, and (2) $df:TU\to TN$ is a vector bundle monomorphism. A related example is the space of immersions $\Imm(-,N):\mathcal{O}(M)\to\Top$, given by $U\mapsto\Imm(U,N)$. This is the space of smooth maps $f:U\to N$ which satisfy (2). We think of an immersion as a local embedding.
\end{example}

The axioms we impose on our functors amount to something like continuity. The first tries to say that our functors should take equivalences to equivalences. At first glace, a category of open subsets of a smooth manifold should have diffeomorphism be the notion of equivalence. Of course, an inclusion map will never be a diffemorphism, so we ask for the next best thing. Let $U,V\in\mathcal{O}(M)$ with $U\subset V$. The inclusion map $i:U\to V$ is called an \emph{isotopy equivalence} if there is an embedding $e:V\to U$ such that the compositions $i\circ e$ and $e\circ i$ are isotopic to the identities of $V$ and $U$ respectively.

\begin{defin}\label{D:goodfunctor}
A contravariant functor $F:\mathcal{O}(M)\to\Top$ is \textsl{good} if
\begin{enumerate}
\item It takes isotopy equivalences to homotopy equivalences, and
\item For any sequence of open sets $U_0\subset U_1\subset\cdots\subset U_i\subset \cdots $, the canonical map $F(\cup_iU_i)\to\holim_{i}F(U_i)$ is a homotopy equivalence.
\end{enumerate}
\end{defin}

Another informal expression of the first axiom is that $F$ behaves well on thickenings. The reader may safely ignore the homotopy limit in the second axiom in favor of this explanation: the functor $F$ is determined by its values on open sets $U$ which are the interior of smooth compact codimension $0$ submanifolds of $M$. Indeed, for any open set $U$, one can select an increasing sequence $V_0\subset V_1\subset\cdots\subset V_i\subset\cdots\subset U$ such that $\cup_iV_i=U$, and each $V_i$ is the interior of a smooth compact codimension $0$ submanifold of $M$. This is a sensible thing to impose in light of our main example of interest, $\Emb(-,N):\mathcal{O}(M)\to\Top$. After all, we are only interested in the values of $\Emb(U,N)$ when $U$ is the interior of some smooth compact manifold. It is also necessary for many of our arguments to assume that $U$ is of this form.

The structure of the category $\mathcal{O}(M)$ is much richer than the usual topology on the real line $\R$, so analogies between functions $f:\R\to\R$ and functors $F:\mathcal{O}(M)\to\Top$ may seem a little weak. Still, there are a few things to say that may be helpful. First, in light of the second axiom above, we could consider the full subcategory of all open sets $U$ which are the interiors of smooth compact codimension $0$ submanifolds of $M$, which we call $\mathcal{O}_{\mathrm{Man}}(M)$. We like to think of $\mathcal{O}_{\mathrm{Man}}(M)\subset\mathcal{O}(M)$ as the analog of the dense subset $\Q\subset\R$ (after all, every continuous function is determined by its values on a dense subset. For more on this, see \refT{charlinear}. We will work almost exclusively in the category $\mathcal{O}_{\mathrm{Man}}(M)$, and so we will make a few remarks about its structure. The objects $U\in\mathcal{O}_{\mathrm{Man}}(M)$ can be coarsely categorized based on their handle dimension. This should be thought of as a more refined notion of dimension of a manifold, and it plays a more important role in this theory than does the ordinary dimension. In particular we will often refer to the handle dimension of an open set $U$, which means the handle dimension of the compact codimension $0$ submanifold whose interior is $U$. Another important subcategory is the full subcategory of open subsets diffeomorphic with at most $k$ open balls. This is the subcategory of $\mathcal{O}_{\mathrm{Man}}(M)$ consisting of those sets $U$ of handle dimension $0$.
\begin{defin}
Let $k\geq 0$. The objects of the full subcategory $\mathcal{O}_k(M)\subset\mathcal{O}(M)$ are those open sets $U$ which are diffeomorphic with at most $k$ open balls in $M$.
\end{defin}
We will return to the categories $\mathcal{O}_k(M)$ later, and their importance will become clear once we define the notion of a polynomial functor.



\section{Derivatives}\label{S:derivatives}

In order to build the Taylor series of a function $f$, we must discuss derivatives. For a smooth function $f:\R\to\R$, its derivative at $0$ is defined by
$$
f'(0)=\lim_{h\to 0}\frac{f(h)-f(0)}{h}.
$$
For our analogy, we will ignore the denominator of the difference quotient in favor of the difference $f(h)-f(0)$. We must decide three things: what plays the role of $0$, what plays the role of $h$, and what plays the role of the difference $f(h)-f(0)$. As for $0$ and $h$, their analogs are, respectively, the empty set $\emptyset$, and the simplest non-empty open set: a set $B$ which is diffeomorphic with an open ball. It is simplest in the sense that it has a handle structure with a single $0$-handle.

As for the difference $f(h)-f(0)$, since $\emptyset\subset B$, for a functor $F$ we have a map $F(B)\to F(\emptyset)$. There are a few ways of computing the difference between two spaces with a map between them. The right thing to do is to compute the homotopy fiber.

\begin{defin}
We define the derivative of $F$ at $\emptyset$ to be
$$
F'(\emptyset)=\hofiber(F(B)\to F(\emptyset)).
$$
\end{defin}

One reason this is natural is because the homotopy fiber, via the long exact sequence in homotopy groups, describes the difference between two spaces in homotopy. If $M$ is connected, then our first axiom (together with a trick allowing us to relate two disjoint open balls in the same path component) implies that the homotopy type of $F'(\emptyset)$ is independent of the choice of $B$.

\begin{example}
Let $F(U)=\Map(U,X)$. Let $B$ be an open ball in $M$. Then $F'(\emptyset)=\hofiber(\Map(B,X)\to\Map(\emptyset,X))\simeq X$, since $\Map(\emptyset,X)=\ast$ and $\Map(B,X)\simeq X$.
\end{example}

\begin{example}
Consider the functor $E(U)=\Emb(U,N)$ and let $B$ be an open ball in $M$. We have $E'(\emptyset)=\hofiber(E(B)\to E(\emptyset))$. An embedding of $B$ is determined by its derivative at a point in $B$ by the inverse function theorem, and so $E(B)$, and hence $E'(\emptyset)$, is equivalent to the space of injective linear maps $\R^m\to \R^n$.
\end{example}

This process can be iterated, just as in ordinary calculus. Choose a basepoint in $F(M)$, which endows $F(U)$ with a basepoint for all $U\in\mathcal{O}(M)$ via the map $F(M)\to F(U)$. For our purposes it is more useful to have formulas for the higher derivatives only in terms of the functor $F$, not its derivatives. Consider the following non-standard formula for the second derivative of $f:\R\to\R$ at $0$:
$$
f''(0)=\lim_{h_1,h_2\to 0}\frac{f(h_1+h_2)-f(h_1)-f(h_2)+f(0)}{h_1h_2}.
$$
Once again for an analogy, will throw away the denominator and focus on the iterated difference
$f(h_1+h_2)-f(h_1)-f(h_2)+f(0)=f(h_1+h_2)-f(h_1)-(f(h_2)-f(0))$. Now all we need is an analog of $+$, for which we will use disjoint union, so $h_1+h_2$ becomes $B_1\coprod B_2$ for two disjoint open balls $B_1, B_2\subset M$. Then we iterate homotopy fibers and define
$$
F''(\emptyset)=\hofiber(\hofiber(F(B_1\coprod B_2)\to F(B_1))\longrightarrow\hofiber(F(B_2)\to F(\emptyset))).
$$
This iterated homotopy fiber is, by definition, the ``total homotopy fiber'' of the following square diagram:
$$\xymatrix{
F(B_1\coprod B_2)\ar[r]\ar[d]&F(B_1)\ar[d]\\
F(B_2)\ar[r]& F(\emptyset)
}
$$
The $k^{th}$ derivative of $F$ at $\emptyset$ is the total homotopy fiber of a $k$-dimensional cubical diagram involving $k$ disjoint open balls. In order to make this precise, we require a brief discussion of cubical diagrams. They are ubiquitous in calculus of functors, and we will use them frequently.

\subsection{Cubical diagrams and total homotopy fibers}\label{S:tfibers}

Details about cubical diagrams can be found in \cite[Section 1]{CalcII}. Other aspects important to this work not appearing in this section have been placed in the appendix to cause minimal distraction. For a finite set $T$, let $|T|$ be its cardinality and $\mathcal{P}(T)$ denote the poset of non-empty subsets of $T$. For instance, if $T=\underline{1}=\{1\}$, this poset looks like $\emptyset \rightarrow \{1\}$, and if $T=\underline{2}$, then we can diagram this poset as a square
$$\xymatrix{
\emptyset\ar[r]\ar[d]& \{1\}\ar[d]\\
\{2\}\ar[r]&\{1,2\}
}
$$
Here we have only indicated those morphisms which are non-identity morphisms and minimal in the sense that they cannot be written as a composition of multiple non-identity morphisms. A $0$-cube is a space, a $1$-cube is a map of spaces and a $2$-cube is a square diagram. In general, the $2^{|T|}$ subsets can be arranged to form a $|T|$-dimensional cube whose edges are the inclusion maps as above. Experience suggests understanding statements for $k$-cubes in the cases $k=2,3$ is usually enough. We will focus almost exclusively on square diagrams.

\begin{defin}
Let $T$ be a finite set. A \textsl{$\abs{T}$-cube} of spaces is a covariant functor $$\mathcal{X}\colon\mathcal{P}(T)\longrightarrow\Top.$$ We may also speak of a cube of based spaces; in this case, the target is $\Top_\ast$.
\end{defin}

We can view a $\abs{T}$-cube $\mathcal{X}$ as a map (i.e.~ a natural transformation of functors) of $(\abs{T}-1)$-cubes $\mathcal{Y}\rightarrow\mathcal{Z}$ as follows. Let $t\in T$. Define $\mathcal{Y}:\mathcal{P}(T-\{t\})\to\Top_\ast$ by $\mathcal{Y}(S)=\mathcal{X}(S)$. Define $\mathcal{Z}:\mathcal{P}(T-\{t\})\to\Top_\ast$ by $\mathcal{Z}(S)=\mathcal{X}(S\cup\{t\})$. There is clearly a natural transformation of functors $\mathcal{Y}\to\mathcal{Z}$, and we may write $\mathcal{X}=(\mathcal{Y}\to\mathcal{Z})$.

\begin{defin}
The \textsl{total homotopy fiber}, or \textsl{total fiber}, of a $\abs{T}$-cube $\mathcal{X}$ of based spaces is the space $\tfiber(\mathcal{X})$ given by the following iterative definition. For a $1$-cube $X_{\emptyset}\to X_{1}$, the total homotopy fiber is defined to be the homotopy fiber of the map $X_{\emptyset}\to X_{1}$. For a $k$-cube $\mathcal{X}$, write it as a map of $(k-1)$-cubes $\mathcal{Y}\to\mathcal{Z}$, and define $\tfiber(\mathcal{X})=\hofiber(\tfiber(\mathcal{Y})\to\tfiber(\mathcal{Z}))$.
\end{defin}

This is well-defined because the homotopy type of $\tfiber(\mathcal{X})$ is independent of the choice of $\mathcal{Y}$ and $\mathcal{Z}$ above by \cite[Proposition 1.2a]{CalcII}. This can be show to be equivalent to the following definition, which is more concise, obviously well-defined, but requires knowledge of homotopy limits.

\begin{prop}\cite[1.1b]{CalcII}
For a $|T|$-cube $\mathcal{X}$ of based spaces, $\tfiber(\mathcal{X})$ is the homotopy fiber of the map 
$$
a(\mathcal{X}):\mathcal{X}(\emptyset)\longrightarrow\underset{S\neq\emptyset}{\holim}\, \, \mathcal{X}(S).
$$ 
\end{prop}

The reader is encouraged to prove this in the case of a square diagram.

\begin{defin}
Let $\mathcal{X}$ be as above.  If $a(\mathcal{X})$ is $k$-connected, we say the cube is \textsl{$k$-cartesian}. In case $k=\infty$, (that is, if the map is a weak equivalence), we say the cube $\mathcal{X}$ is \textsl{homotopy cartesian}.
\end{defin}

For a space ($0$-cube) $X$, the convention is that $k$-cartesian means $(k-1)$-connected, and  for a map ($1$-cube) $X\to Y$ to be $k$-cartesian means it is $k$-connected, so its homotopy fibers are $(k-1)$-connected. A square 
$$\xymatrix{
X_\emptyset\ar[r]\ar[d]& X_1\ar[d]\\
X_2\ar[r] & X_{12}
}
$$
is homotopy cartesian if the map $X_\emptyset\to\holim(X_1\to X_{12}\leftarrow X_2)$ is a homotopy equivalence. Such a square is often referred to as a \emph{homotopy pullback} square because $\holim(X_1\to X_{12}\leftarrow X_2)$ is the space of all $(x_1,\gamma, x_2)$ such that $x_i\in X_i$ for $i=1,2$ and $\gamma$ is a path in $X_{12}$ between the images of $x_1$ and $x_2$. In contrast, the pullback of $X_1\to X_{12}\leftarrow X_2$ is the space of all $(x_1,x_2)\in X_1\times X_2$ such that the images of the $x_i$ in $X_{12}$ are equal. There is a useful relationship between pullbacks and homotopy pullbacks. If
$$\xymatrix{
X_\emptyset\ar[r]\ar[d]& X_1\ar[d]\\
X_2\ar[r] & X_{12}
}
$$
is a pullback square, then it is a homotopy pullback if either $X_1\to X_{12}$ or $X_2\to X_{12}$ is a fibration. That is, in this case the map from the pullback to the homotopy pullback is an equivalence. A similar criterion can be formulated for general cubes, though it is more complicated. A useful and familiar example of a homotopy pullback is obtained by setting $X_2$ to a point and letting $X_1\to X_{12}$ be a fibration whose fiber over the image of $X_2$ in $X_{12}$ is $X_\emptyset$.

Viewing a $(|T|+1)$-cube $\mathcal{Z}$ as a map of $|T|$-cubes $\mathcal{X}\to\mathcal{Y}$ as in our iterative definition of total homotopy fiber, choose a basepoint $y\in\mathcal{Y}(\emptyset)$, which bases each $\mathcal{Y}(S)$, and define a $|T|$-cube $\mathcal{F}_y(S)=\hofiber(\mathcal{X}(S)\to\mathcal{Y}(S))$.

\begin{prop}\label{P:cartesianfibercube}\cite[1.18]{CalcII}
With $\mathcal{X}, \mathcal{Y}, \mathcal{Z}$ as above, the $(|T|+1)$-cube $\mathcal{X}$ is $k$-cartesian if and only if for each choice of basepoint $y\in\mathcal{Y}(\emptyset)$, the $|T|$-cube $S\mapsto\mathcal{F}_y(S)$ is $k$-cartesian.
\end{prop}

For $|T|=1$, this says that a map of spaces $X\to Y$ is $k$-connected if and only if all of its homotopy fibers are $k$-cartesian, which means $(k-1)$-connected in the case of a $0$-cube. We present one final fact which will be useful in the proof of \refT{etotone}.

\begin{prop}\label{P:ScubeTcartesian}\cite[1.22]{CalcII}
Let $\mathcal{X}, \mathcal{Y}$ be $|T|$-cubes, and suppose we have a map $\mathcal{X}\to\mathcal{Y}$ such that for all $S\neq\emptyset$, $\mathcal{X}(S)\to\mathcal{Y}(S)$ is $k$-connected. Then the map $\holim_{S\neq\emptyset}\mathcal{X}(S)\to\holim_{S\neq\emptyset}\mathcal{Y}(S)$ is $(k-|T|+1)$-connected.
\end{prop}


Returning to our discussion of derivatives, we can now make a sensible definition of the derivatives of $F$ at $\emptyset$.

\begin{defin} 
Let $B_1,\ldots, B_k$ be pairwise disjoint open balls in $M$. Define a $k$-cube of spaces by the rule $S\mapsto F(\cup_{i\notin S}B_i)$. Define the \textsl{$k^{th}$ derivative of $F$ at the empty set}, denoted $F^{(k)}(\emptyset)$, to be the total homotopy fiber of the $k$-cube $S\mapsto F(\cup_{i\notin S}B_i)$.
\end{defin}

\begin{example}
We can compute the derivatives of $F(U)=\Map(U,X)$. We have already seen that $F'(\emptyset)\simeq X$. Let $B_1,B_2$ be disjoint open balls. $F''(\emptyset)$ is the total homotopy fiber of the square
$$\xymatrix{
F(B_1\coprod B_2)\ar[r]\ar[d]& F(B_1)\ar[d]\\
F(B_2)\ar[r] & F(\emptyset)
}
$$
Since each $B_i$ is homotopy equivalent to a point $\ast_i$, and $\Map(-,X)$ preserves homotopy equivalences, this is equivalent to the total homotopy fiber of the square
$$\xymatrix{
F(\ast_1\coprod \ast_2)\ar[r]\ar[d]& F(\ast_1)\ar[d]\\
F(\ast_2)\ar[r] & F(\emptyset)
}
$$
Clearly $F(\ast_1\coprod \ast_2)=X\times X$, and by our calculation above, we see that $F''(\emptyset)$ is the total homotopy fiber of the square
$$\xymatrix{
X\times X\ar[r]\ar[d]& X\ar[d]\\
X\ar[r] & \ast
}
$$
Here the vertical map $X\times X\to X$ is projection onto the second coordinate, and the horizontal map is projection onto the first coordinate. Using our iterative definition of homotopy fiber (and taking fibers vertically), we see that $F''(\emptyset)=\hofiber(X\stackrel{id}{\to} X)\simeq\ast$. Alternately, we could observe that this square is both a pullback and a homotopy pullback. A similar computation shows that $F^{(k)}(\emptyset)\simeq\ast$ for $k\geq 3$. That is, all derivatives but the first of $F(U)=\Map(U,X)$ vanish, which suggests this should be a linear functor. It is, as we will see in the next section.
\end{example}

\begin{example}
Let us compute the first two derivatives of $F(U)=\Map(U^2,X)$. We have $F'(\emptyset)=\hofiber(\Map(B^2,X)\to\Map(\emptyset,X))\simeq \hofiber(X\to\ast)\simeq X$. $F''(\emptyset)$ is the total homotopy fiber of the square
$$\xymatrix{
\Map((B_1\coprod B_2)^2,X)\ar[r]\ar[d]& \Map(B_1^2,X)\ar[d]\\
\Map(B_2^2,X)\ar[r] & \Map(\emptyset,X)
}.
$$
Since $\Map((B_1\coprod B_1)^2,X)=\Map(B_1^2,X)\times\Map(B_2^2,X)\times \Map(B_1\times B_2,X)\times \Map(B_2\times B_1,X)$, we have $F''(\emptyset)\simeq\Map(B_1\times B_2,X)\times\Map(B_2\times B_1,X)\simeq  X^2$. All of the higher derivatives are contractible. In a similar fashion, one can compute the first $k$ derivatives of $F(U)=\Map(U^k,X)$; all derivatives of order greater than $k$ are contractible.
\end{example}

\subsection{Criticism of analogies}

We justify our definition of derivative based on the classification theorem for homogeneous functors which appears below as \refT{homogclass}, in which the derivatives at the empty set play a central role. Despite the importance of the derivatives as we have defined them, we have reason to be interested in the derivative of $F$ at an arbitrary open set. We even have reason to be interested in something which formally resembles a derivative (the homotopy fiber of a restriction map) as described above, but which does not simply involve studying differences based on taking disjoint unions with open balls. It is natural to make the following definition.

\begin{defin}
For an open set $V\in\mathcal{O}(M)$ and an open ball $B$ disjoint from $V$, define
$$
F'(V)=\hofiber(F(V\coprod B)\to F(V)).
$$
\end{defin}

Although the disjoint union is an obvious candidate for the analog of sum, it is not at all clear that we should ignore more general unions, for example, the attaching of a handle. In fact, we should not. As we have mentioned, it is enough for us to understand the values of a functor on open sets $V$ which are the interior of a compact codimension $0$ submanifold $L$ of $M$; that is, when $V\in\mathcal{O}_{\mathrm{Man}}(M)$. For the purposes of this informal discussion, we will replace $V$ with $L$. We wish to consider more generally $\hofiber(F(L\cup_f H^i)\to F(L))$. The special case of $i=0$ is the disjoint union of $L$ with an $m$-dimensional disk. Similar criticisms apply to the study higher derivatives. More general differences will become important when we tackle the question of convergence and the analog of a bound on $f^{(k+1)}(x)$ for $x$ close to zero, which is important in understanding the remainder $R_k(x)=|f(x)-T_kf(x)|$.

\section{Polynomial Functors}\label{S:polynomials}

A function $f:\R\to\R$ is \emph{linear} if $f(x+y)=f(x)+f(y)$ for all $x,y$. More generally, we might say a function is linear if $f(x+y)-f(x)-f(y)+f(0)=0$. Making analogies as we did in \refS{derivatives}, and being more flexible about the analog of sum (and using an arbitrary union in place of the disjoint union), this leads one to say that a functor $F:\mathcal{O}(M)\to\Top$ is linear if for all open $V,W$ in $M$ the total homotopy fiber of 
$$\xymatrix{
F(V\cup W)\ar[r]\ar[d] & F(W)\ar[d]\\
F(V)\ar[r]&F(V\cap W)
}
$$
is contractible. This implies that the second (and higher) derivatives of $F$ vanish by letting $V$ and $W$ be disjoint open balls, but linearity is clearly a stronger condition. Linear functors are also called \emph{polynomial of degree $\leq 1$}, or \emph{excisive}. We pause for an example before formalizing this definition.

\begin{example}
Let $X$ be a space. The functor $U\mapsto\Map(U,X)$ is linear. This follows from the fact that $\Map(-,X)$ sends (homotopy) pushout squares to (homotopy) pullback squares. See \refP{maphocolimtoholim}.
\end{example}

We can reformulate this in a way more suitable to our needs, and although it may seem a bit strange at first, the proof of \refT{charlinear} should help the reader understand why the definition is presented this way.
\begin{defin}\label{D:linear}
A functor $F:\mathcal{O}(M)\to\Top$ is \textsl{polynomial of degree $\leq1$} if for all $U\in\mathcal{O}(M)$ and for all disjoint nonempty closed subsets $A_0,A_1\subset U$, the diagram
$$\xymatrix{
F(U)\ar[r]\ar[d] & F(U-A_0)\ar[d]\\
F(U-A_1)\ar[r]&F(U-(A_0\cup A_1))
}
$$
is homotopy cartesian.
\end{defin}

To relate this back to the definition above, note that if we put $W=U-A_0$ and $V=U-A_1$, then $U=V\cup W$, and $U-(A_0\cup A_1)=V\cap W$. The reason for this is that it is convenient for the purposes of inductive arguments (we will see this first in the proof of \refT{charlinear}) to think about ``punching holes'' in an open set to reduce its handle dimension. The definition of polynomial of higher degree generalizes the notion of linearity.

\begin{defin}\label{D:polynomial}
A functor $F:\mathcal{O}(M)\to\Top$ is called \textsl{polynomial of degree $\leq k$} if for all $V\in\mathcal{O}(M)$ and for all pairwise disjoint nonempty closed subsets $A_0,A_1,\ldots, A_{k+1}\subset V$, the map
$F(V)\to\holim_{S\neq\emptyset}F(V-\cup_{i\in S}A_i)$
is a homotopy equivalence; equivalently, the diagram $S\mapsto F(V-\cup_{i\in S}A_i)$ is homotopy cartesian.
\end{defin}

To compare this with the definition of the $k$th derivative, let $V$ be $k+1$ disjoint open balls and let the $A_i$ be the components of $V$. Thus a polynomial of degree $\leq k$ has contractible derivatives of order $k+1$ and above. 

\begin{prop}
If $F$ is polynomial of degree $\leq k$, then it is polynomial of degree $\leq k+1$.
\end{prop}

This is certainly something that had better be true if this definition is to make any sense. It is not completely trivial, but follows from the fact that if two opposing $(k+1)$-dimensional faces of a $(k+2)$-cube are homotopy cartesian, then that $(k+2)$-cube is itself homotopy cartesian. Now let us consider several more examples.

\begin{example}
The functor $U\mapsto \Map(U^k,X)$ is polynomial of degree $\leq k$ (but not polynomial of lower degree). This basically follows from the pigeonhole principle. Let $A_0,\ldots, A_k$ be pairwise disjoint nonempty closed subsets of $U$. For a point $(x_1,\ldots, x_k)\in U^k$, each $x_i$ is in at most one $A_j$, hence there is some $l$ such that $x_i\in U-A_l$ for all $i$ by the pigeonhole principle. Therefore $U^k=\cup_{i=1}^k(U-A_i)^k$. It follows immediately that, $U^k=\colim_{S\neq\emptyset}(U-\cup_{i\in S}A_i)^k$, and one can show that in fact $U^k\simeq\hocolim_{S\neq\emptyset}(U-\cup_{i\in S}A_i)^k$. Since $\Map(-,X)$ preserves equivalences and turns homotopy colimits into homotopy limits (\refP{maphocolimtoholim}), we have an equivalence $\Map(U^k,X)\simeq\holim_{S\neq\emptyset}\Map((U-\cup_{i\in S}A_i)^k,X)$.
\end{example}

\begin{example}
We can generalize the previous example without doing any extra work as 
follows. Let $C(k,U)\subset U^k$ be the configuration space of $k$ points 
in $U$ (those $(x_1,\ldots, x_k)\in U^k$ such that $x_i\neq x_j$ for 
$i\neq j$). The group $\Sigma_k$ acts on $C(k,U)$ by permuting the 
coordinates, and we let $\binom{U}{k}=C(k,U)/\Sigma_k$ denote the 
quotient by this action. This gives us the space of unordered 
configurations of $k$ points in $U$. The same argument as in 
the previous example shows that both $U\mapsto\Map(C(k,U),X)$ 
and $U\mapsto\Map(\binom{U}{k},X)$ are polynomial of degree $\leq k$.
\end{example}

\begin{example}
The functor $U\mapsto \Emb(U,N)$ is not polynomial of degree $\leq k$ for any $k$. We will indicate why for $k=1$. Let $A_0,A_1\subset U$ be pairwise disjoint closed subsets, and put $U_i=U-A_i$, and $U_{12}=U_1\cap U_2$. We are asked to check whether the map $\Emb(U_1\cup U_2,N)\to\holim(\Emb(U_1,N)\to\Emb(U_{12},N)\leftarrow\Emb(U_2,N))$ is an equivalence. That is, given $f_i\in\Emb(U_i,N)$ with a homotopy between their restrictions to $U_{12}$, is this enough to determine an element of $\Emb(U_1\cup U_1,N)$? It is not, due to an obstruction, namely that $f_1(U_1)$ and $f_2(U_2)$ might intersect in $N$. It is, however, true that the map $\Emb(U_1\cup U_2,N)\to\holim(\Emb(U_1,N)\to\Emb(U_{12},N)\leftarrow\Emb(U_2,N))$ has a certain connectivity; see \refS{disj} and \refT{tonedisj}
\end{example}

\begin{example}
The functor $U\mapsto\Imm(U,N)$ is a polynomial of degree $\leq 1$. Let $A_0,A_1\subset U$ be pairwise disjoint closed, and put $U_i=U-A_i$, and $U_{01}=U_0\cap U_1$. Then the square
$$\xymatrix{
\Imm(U_0\cup U_1,N)\ar[d]\ar[r] & \Imm(U_0,N)\ar[d]\\
\Imm(U_1,N)\ar[r] & \Imm(U_{01},N)\\
}
$$
is clearly a pullback, since being an immersion is a local condition, and immersions of $U_0$ and $U_1$ which agree on their intersection make an immersion of the union. It is a homotopy pullback because the restriction map $\Imm(U_0,N)\to \Imm(U_{01},N)$ is a fibration. This fact is a reformulation of the Smale-Hirsch theorem. This isn't quite technically correct; the Smale-Hirsch theorem does not apply to the restriction map of open sets. However, this can be overcome without too much difficulty. See \cite[Lemma 1.5]{W:EI1}.
\end{example}

\subsection{Characterization of polynomials}\label{S:polyclass}

\refT{polychar} below is a structure theorem for polynomials, and later we will discuss a structure theorem for homogeneous polynomials, \refT{homogclass}. \refT{charlinear}, a structure theorem for linear functors (polynomials of degree $\leq 1$), which contains aspects of the proofs of both \refT{polychar} and \refT{homogclass}, will be given below, and it has a simple parallel for ordinary linear functions $f:\R\to \R$. The techniques of its proof are used many times in this paper.

Consider the following proof that every continuous linear function $f:\R\to\R$ is of the form $f(x)=ax$. Let $a=f(1)$. Linearity implies $f(n)=an$ for $n$ a natural number. If $p$ and $q$ are natural numbers with $q\neq 0$, then $ap=f(q\frac{p}{q})=qf(\frac{p}{q})$ by linearity, and so $f(\frac{p}{q})=a\frac{p}{q}$. By density of $\Q$ in $\R$ and continuity of $f$, this implies $f(x)=ax$ for all real $x$.

Let $p:Z\to M$ be a fibration, and let $\Gamma(M,Z;p)$ be its space of sections. For example, if $Z=M\times X$ and $p$ is the projection, $\Gamma(M,Z;p)=\Map(M,X)$. The following theorem says that all linear functors $F$ such that $F(\emptyset)=\ast$ are the space of sections of some fibration. Or, more roughly, that they are all (twisted) mapping spaces.

\begin{thm}\label{T:charlinear}
Let $F:\mathcal{O}(M)\to\Top$ be a good functor such that $F(\emptyset)=\ast$ and which is polynomial of degree $\leq1$. Then there is a fibration $p:Z\to M$ for some space $Z$ and a natural transformation $F(U)\to\Gamma(U,Z;p)$ which is an equivalence for all $U\in\mathcal{O}(M)$.
\end{thm}
\begin{proof}
First we will make the natural transformation $F(U)\to\Gamma(U,Z;p)$. Let $\mathcal{O}^{(1)}(V)$ denote the category of open subsets of $V$ which are diffeomorphic to exactly one open ball. Note that all inclusions in this category are isotopy equivalences, and that the realization $|\mathcal{O}^{(1)}(V)|\simeq V$. Let $Z=\hocolim_{U\in\mathcal{O}^{(1)}(M)}F(U)$. Since $F$ takes isotopy equivalences to homotopy equivalences, $Z$ quasifibers over $|\mathcal{O}^{(1)}(M)|\simeq M$ with space of sections equivalent to $\holim_{U\in\mathcal{O}^{(1)}(M)}F(U)$ by  \refT{holimquasifiber}. There is a natural transformation $F(V)\to\holim_{U\in\mathcal{O}^{(1)}(V)}F(U)$ since $F(V)\simeq\holim_{U\in\mathcal{O}(V)}F(U)$ by \refT{holiminitial} and $\mathcal{O}^{(1)}(M)\to\mathcal{O}(M)$ induces the map in question. We define $\Gamma(V)=\holim_{U\in\mathcal{O}^{(1)}(V)}F(U)$. We now must show $F(V)\to\Gamma(V)$ is an equivalence. To do so, it is enough by the second part of \refD{goodfunctor} to check that it is an equivalence when $V$ is the interior of a compact codimension zero submanifold $L$ of $M$. We will proceed by induction on the handle dimension of $V$.

Let $k$ be the handle dimension of $L$. The base case to consider is $k=0$, when $V$ is a disjoint union of finitely many open balls. For this, we will induct on the number of components $l$ of $V$. The base case is $l=1$, and in this case $V$ is a final object in the category $\mathcal{O}^{(1)}(V)$, and so the map $F(V)\to\Gamma(V)$ is an equivalence by \refT{holiminitial}. For $l>1$, let $A_0,A_1$ be two distinct components of $V$, and put $V_S=V-\cup_{i\in S}A_i$ for $S\subset \{0,1\}$. Consider the following diagram
$$\xymatrix{
F(V)\ar[r]\ar[d]& \Gamma(V)\ar[d]\\
\holim_{S\neq \emptyset}F(V_S)\ar[r] & \holim_{S\neq \emptyset}\Gamma(V_S)
}
$$
Since both $F$ and $\Gamma$ are polynomial of degree $\leq 1$, the vertical maps are equivalences, and by induction, each map $F(V_S)\to\Gamma(V_S)$ is an equivalence, and hence the induced map of homotopy limits over $S$ is an equivalence by \refT{holimheq}. Therefore the top arrow is an equivalence as well.

The general case proceeds in a similar fashion. Let $k>0$ be the handle dimension of $V$, and let $l$ denote the number of handles of dimension $k$. Let $e:D^k\times D^{m-k}\to L$ be one of these $k$-handles. Let $D_0,D_1\subset D^k$ be disjoint disks, and put $A'_i=D_i\times D^{m-k}$. Then $A_0=V\cap A'_0$ and $A_1=V\cap A'_1$ are nonempty disjoint closed subsets of $V$, and if we put $V_S=V-\cup_{i\in S}A_i$, then for $S\neq\emptyset$, $V_S$ is the interior of a compact codimension zero submanifold $L_S$ which can be given a handle structure with fewer than $l$ handles of dimension $k$ (see \ref{handles} for a picture in a slightly different case). Once again consider the following diagram.
$$\xymatrix{
F(V)\ar[r]\ar[d]& \Gamma(V)\ar[d]\\
\holim_{S\neq \emptyset}F(V_S)\ar[r] & \holim_{S\neq \emptyset}\Gamma(V_S)
}
$$
The vertical arrows are equivalences because $F$ and $\Gamma$ are polynomial of degree $\leq 1$. For $S\neq\emptyset$, the map $F(V_S)\to \Gamma(V_S)$ is an equivalence by induction on $l$, and hence so is the bottom horizontal arrow. It follows that the top arrow is an equivalence as well.
\end{proof}

\begin{rem}
The idea of this proof is philosophically similar to that which classifies continuous linear functions. We first constructed the desired functor $\Gamma$ by averaging (taking a homotopy limit) the values of $F$ on single open balls (akin to $a=f(1)$; we took an average to ensure functoriality), and we see from the proof that $\Gamma$, and hence $F$, is completely determined by the value of $F$ on an open ball. Then we showed using linearity with a handle induction argument that this implied that $F(V)\to\Gamma(V)$ was an equivalence for $V\in\mathcal{O}_{\mathrm{Man}}(M)$ (our analog of $\Q$). Finally we used continuity to conclude the result for general open sets $V$.
\end{rem}

\begin{rem}
We have already seen that $\Imm(M,N)$ is polynomial of degree $\leq 1$. We may ask how to express it as a space of sections. In this case, an immersion $f$ is a section of a bundle over $M$ whose fiber at $x\in M$ is the space of vector bundle monomorphisms $T_xM\to T_{f(x)}N$. This is, once again, a version of the Smale-Hirsch Theorem.
\end{rem}

A proof similar to that in \refT{charlinear} characterizes polynomials in terms of their values on finitely many open balls, and it also utilizes a similar handle induction argument.

\begin{thm}\cite[Theorem 5.1]{W:EI1}\label{T:polychar}
Suppose $F_1\to F_2$ is a natural transformation of good functors and that $F_i$ is a polynomial of degree $\leq k$ for $i=1,2$. If $F_1(V)\to F_2(V)$ is an equivalence for all $V\in\mathcal{O}_k(M)$, then it is an equivalence for all $V\in \mathcal{O}(M)$.
\end{thm}

Note that a polynomial $p:\R\to\R$ of degree $k$ such that $p(0)=0$ is determined by its values on $k$ distinct points; similarly, our polynomial functors $F$ are completely determined by their values on the category of at most $k$ open balls. 

\subsection{Approximation by polynomials}\label{S:polyapprox}

Now we will construct the $k$th Taylor polynomial $T_kF$ for a functor $F$. Proceeding with an ordinary Taylor polynomial in mind, we would like to construct a functor $T_kF$ which has the following properties:
\begin{itemize}
\item The derivatives $F^{(i)}(\emptyset)$ and $(T_kF)^{(i)}(\emptyset)$ agree for $0\leq i\leq k$.
\item $T_kF$ is polynomial of degree $\leq k$.
\item There is a natural transformation $F\to T_kF$, so that we may discuss the ``remainder'' $R_kF=\hofiber(F\to T_kF)$.
\end{itemize}

Looking back at our discussion of derivatives, we computed $F^{(i)}(\emptyset)$ by looking at the total homotopy fiber of a cubical diagram of the values of $F$ on at most $i$ disjoint open balls. One way to ensure that the derivatives of order at most $k$ of $F$ and $T_kF$ agree is to make the values of $F(V)$ and $T_kF(V)$ agree when $V$ is a disjoint union of at most $k$ open balls. With this in mind, for $V\in\mathcal{O}(M)$, recall the poset $\mathcal{O}_k(V)$ of open subsets of $U$ which are diffeomorphic with at most $k$ open balls in $V$. It is a subposet of $\mathcal{O}(V)$, and we want the values of $F$ and $T_kF$ to agree on these subcategories. 

\begin{defin}\label{D:tkf}
Let
$$T_kF(V)=\holim_{U\in\mathcal{O}_k(V)}F(U).$$
\end{defin}

This is a (homotopy) \emph{Kan extension} of $F$ along the inclusion of the subcategory $\mathcal{O}_k(V)\to\mathcal{O}(V)$. It says that the value of $T_kF$ at a given open set $V$ is an ``average'' of the values of $F$ on at most $k$ open balls contained in $V$. Note that if $V$ itself is diffeomorphic with at most $k$ open balls, then $V$ is a final object in $\mathcal{O}_k(V)$, and so $T_kF(V)=\holim_{U\in\mathcal{O}_k(V)}F(U)\simeq F(V)$, so we really have correctly prescribed the values of $T_kF$ the way we said we would. 

It is not clear from \refD{tkf} that $T_kF$ is a polynomial of degree $\leq k$, but it turns out that this is so. The proof is not trivial. Let us content ourselves with knowledge  that an ordinary polynomial of degree $k$ such that $p(0)=0$ is completely determined by its values on at most $k$ points, and it is clear from the definition of $T_kF$ as an extension over the subcategory of at most $k$ ``points'' that the analog of this is true.

There is a natural transformation $F\to T_kF$ given by observing that the inclusion $\mathcal{O}_k(V)\to\mathcal{O}(V)$ induces a map of homotopy limits
$$
F(V)\simeq\holim_{U\in\mathcal{O}(V)}F(U)\to \holim_{U\in\mathcal{O}_k(V)}F(U)=T_kF(V)
$$
and noting that the first equivalence follows since $V$ is a final object in $\mathcal{O}(V)$ (see \refT{holiminitial} in the Appendix).

Note that $\mathcal{O}_0(V)$ contains only the empty set for all $V$, and so $T_0F(V)=F(\emptyset)$ for all $V$.

\begin{example}
Since $F(V)=\Map(V,X)$ is polynomial of degree $\leq 1$, $F(V)\to T_1F(V)$ is an equivalence by \refT{polychar}, since their values agree when $V$ is a single open ball.
\end{example}

\begin{example}
The linearization of embeddings is immersions. That is, $T_1\Emb(V,N)\simeq\Imm(V,N)$. The natural transformation $\Emb(V,N)\to\Imm(V,N)$ is an equivalence when $V$ is a single open ball, and hence $T_1\Emb(U,N)=\holim_{V\in\mathcal{O}_1(U)}\Emb(V,N)$ is equivalent to $\holim_{V\in\mathcal{O}_1(U)}\Imm(V,N)\simeq\Imm(U,N)$, with the last equivalence given by the fact that $\Imm(-,N)$ is polynomial of degree $\leq 1$, as in the previous example.
\end{example}

\subsection{The Taylor Tower}\label{S:tower}

Armed with a definition of $T_kF$, we can now form the ``Taylor tower'' of $F$, the analog of the Taylor series. The inclusion $\mathcal{O}_{k-1}(V)\to \mathcal{O}_{k}(V)$ induces a map $T_kF(V)\to T_{k-1}F(V)$, and so we obtain a tower of functors
$$
\cdots\to T_kF\to T_{k-1}F\to\cdots\to T_1F\to T_0F.
$$
Since $V$ is a final object in $\mathcal{O}(V)$, we may identify $F(V)$ with $\holim_{\mathcal{O}(V)}F$, and the inclusion $\mathcal{O}_k(V)\to\mathcal{O}(V)$ induces maps $F\to T_kF$ which are compatible with one another. Hence there is a natural transformation $F\to\holim_k T_kF$, and we would like to know under what circumstances this map is an equivalence; that is, when the Taylor tower of $F$ converges to $F$. This is the subject of \refS{conv}. Before we embark on questions of convergence, it will be useful to understand the differences $\hofiber(T_kF\to T_{k-1}F)$.

\section{Homogeneous Functors}\label{S:homogclass}

An explicit description of polynomial functors is perhaps too much to hope for, so we will content ourselves with a classification of homogenous functors. Fortunately there is a parallel with ordinary calculus here too. For $f:\R\to\R$, consider the $k$th homogeneous piece of its Taylor series, $L_kf(x)=T_kf(x)-T_{k-1}f(x)=f^{(k)}(0)\frac{x^k}{k!}$. The classification of homogeneous functors shares a similar form. Roughly speaking, it is the space of sections of a fibration over $\binom{M}{k}$ whose fibers are the derivatives $F^{(k)}(\emptyset)$. We will state this more precisely below, but first we define what it means for a functor to be homogeneous and consider some examples.

\begin{defin}
A functor $E:\mathcal{O}(M)\to\Top$ is \textsl{homogeneous of degree $k$} if it is polynomial of degree $\leq k$ and $T_{k-1}E(V)\simeq\ast$ for all $V$.
\end{defin}

\begin{example}
For a good functor $F$, choose a basepoint in $F(M)$. This bases $F(V)$ for all $V\in\mathcal{O}(M)$. The functor $L_kF=\hofiber(T_kF\to T_{k-1}F)$ is homogeneous of degree $k$. That it is polynomial of degree $\leq k$ follows from the fact that $T_kF$ and $T_{k-1}F$ are both polynomial of degree $\leq k$. To see that $T_{k-1}L_kF(V)\simeq\ast$ for all $V$, first observe that $T_{k-1}$ commutes with homotopy fibers (see \refT{holimcommute}; homotopy limits commute), and next observe that $T_{k-1}T_kF\simeq T_{k-1}F$. Indeed, $T_{k-1}T_kF(V)=\holim_{W\in\mathcal{O}_{k-1}(V)}\holim_{U\in\mathcal{O}_{k}(W)}F(U)$, and since $W$ is diffeomorphic with at most $k-1$ open balls, it is a final object in $\mathcal{O}_{k}(W)$, and so $\holim_{U\in\mathcal{O}_{k}(W)}F(U)\simeq F(W)$.
\end{example}

\begin{example}
The functor $U\mapsto\Map(U^2,X)$ is polynomial of degree $\leq 2$, so its quadratic approximation $T_2\Map(U^2,X)\simeq\Map(U^2,X)$. However, it is not homogeneous of degree $2$, because, as we showed above, it has a non-trivial first derivative, which would necessarily vanish were it homogeneous. In fact, $T_1\Map(U^2,X)\simeq\Map(U,X)$. Let $U\to U^2$ be the diagonal map. This gives rise to a restriction $\Map(U^2,X)\to\Map(U,X)$. Note that when $U$ is a single open ball, $\Map(U^2,X)\to\Map(U,X)$ is an equivalence, and since $\Map(U,X)$ is polynomial of degree $\leq 1$, it follows from \refT{polychar} that $T_1\Map(U^2,X)\simeq\Map(U,X)$. Therefore $L_2\Map(U^2,X)=\hofiber(\Map(U^2,X)\to\Map(U,X))$. Similarly, $U\mapsto\Map(U^k,X)$ is not homogeneous of degree $k$ unless $k=1$.
\end{example}

\begin{example}
Let us compute $L_3\Map(U^3,X)=\hofiber(T_3\Map(U^3,X)\to T_2\Map(U^3,X))$. As in the previous example, $T_3\Map(U^3,X)\simeq\Map(U^3,X)$. Let $\Delta(U)\subset U^3$ denote the fat diagonal. $\Delta(U)=\{(x_1,x_2,x_3)| x_i=x_j\mbox{ for some }i\neq j\}$. We would like to claim that $U\mapsto\Map(\Delta(U),X)$ is a model for $T_2\Map(U^3,X)$, and while this is in spirit the case, our answer will be slightly different. 

We proceed as follows: For $S\subset\{1,2,3\}$, let $\Delta_S(U)=\{(x_1,x_2,x_3) | x_i=x_j\mbox{ for all }i,j\in S\}$. Then $\Delta(U)=\colim_{1<|S|}\Delta_S(U)$ (the union of these spaces covers $\Delta(U)$, and we define $\tilde{\Delta}(U)=\hocolim_{1<|S|}\Delta_S(U)$. Thus, since $\Map(-,X)$ turns homotopy colimits into homotopy limits by \refP{maphocolimtoholim}, hence $\Map(\tilde{\Delta}(U),X)=\holim_{1<|S|}\Map(\Delta_S(U),X)$. It is clear that $\Map(\Delta_S(U),X)$ is a polynomial of degree $\leq 4-|S|$, and since $1<|S|\leq 3$, for all $S$ under consideration, $\holim_{1<|S|}\Map(\Delta_S(U),X)$ is polynomial of degree $\leq 2$ because each functor in the diagram is polynomial of degree $\leq 2$. Note that $U^3\simeq\hocolim_{1\leq |S|}\Delta_S(U)$, and hence there is a natural transformation of functors $\Map(U^3,X)\to\Map(\tilde{\Delta}(U),X)$ given by the obvious inclusion of categories. 

By inspection, when $U$ is a union of at most two open balls, the map $\Map(U^3,X)\to\Map(\tilde{\Delta}(U),X)$ is an equivalence, and so by \refT{polychar}, $T_2\Map(U^3,X)\simeq\Map(\tilde{\Delta}(U),X)$. It follows that $L_3\Map(U^3,X)\simeq\hofiber(\Map(U^3,X)\to\Map(\tilde{\Delta}(U),X))$.
\end{example}

Spaces of maps are special cases of sections of bundles, and we can generalize further to include examples such as these.

\begin{example}
Let $p:Z\to\binom{M}{k}$ be a fibration with a section. 
Let $\Gamma(\binom{M}{k},Z;p)$ denote its (based) space of sections. The assignment $U\mapsto\Gamma(\binom{U}{k},Z;p)$ is polynomial of degree $\leq k$. Define
$$\Gamma\left(\del\binom{U}{k},Z;p\right)=\underset{N\in\mathcal{N}}{\hocolim}\;\Gamma\left(\binom{U}{k}\cap Q,Z;p\right).$$
One may think of this as the space of germs of sections near the fat diagonal. It turns out that $T_{k-1}\Gamma(\binom{U}{k},Z;p)\simeq\Gamma(\del\binom{U}{k},Z;p)$, and hence 
$$
\Gamma^c\left(\binom{U}{k},Z;p\right)=\hofiber\left(\Gamma\left(\binom{U}{k},Z;p\right)\rightarrow\Gamma\left(\del\binom{U}{k},Z;p\right)\right).
$$
is homogeneous of degree $k$. We refer to $\Gamma^c$ as the space of compactly supported sections.
\end{example}

\subsection{Classification of homogeneous polynomials}

The last example in the previous section is quite general, according to the classification of homogeneous functors. 




\begin{thm}[\cite{W:EI1}, Theorem 8.5]\label{T:homogclass}\ \ 
\item Let $E$ be homogeneous of degree $k$. Then there is an equivalence, natural in $U$, 
$$E(U)\longrightarrow \Gamma^c\left(\binom{U}{k},Z;p\right),$$ 
where $\Gamma^c$ is the space of compactly supported sections of a fibration $p:Z\rightarrow\binom{U}{k}$. The fiber over $S$ of the fibration $p$ is the total homotopy fiber of a $k$-cube of spaces made up of the values of $E$ on a tubular neighborhood of $S$. In particular, if $E(U)=\hofiber(T_kF(U)\to T_{k-1}F(U))$, then the fibers of the classifying fibration are the derivatives $F^{(k)}(\emptyset)$.
\end{thm}

This has a pleasing analogy with the $k$th homogeneous piece $\frac{x^k}{k!}f^{(k)}(0)$ of the Taylor series centered at $0$ for a smooth function $f$, 
where $\binom{U}{k}$ plays the role of $\frac{x^k}{k!}$, and, of course, $F^{(k)}(\emptyset)$ plays the role of $f^{(k)}(0)$. We will not discuss the proof of \refT{homogclass}, but remark that most of the required tools are on display in the proof of \refT{charlinear}. The classifying fibration $p:Z\to\binom{U}{k}$ is the pullback of a fibration $p:Z\to\binom{M}{k}$, induced by the inclusion $U\to M$.

\section{Convergence and Analyticity}\label{S:conv}

Now that we can construct a Taylor tower for a functor $F$ and understand a bit about its structure, we are ready to ask whether or not it approximates the functor $F$ in a useful way. The Taylor series of a function $f:\R\to\R$ need not converge to $f$; in fact, the series need not converge at all. We will discuss the extent to which an approximation by polynomial functors does a suitable job approximating the homotopy type of the values of a given functor. The reader may already suspect that a ``suitable'' approximation is one which approximates the homotopy type of through a range. On $\R$, $|x-y|$ measures the difference of $x$ and $y$, and in $\Top$, a useful ``metric'' for measuring the difference between spaces $X$ and $Y$ with respect to a map $f:X\to Y$ is to ask for the connectivity of the homotopy fiber $\hofiber(f)$.

Two natural questions to ask are:
\begin{enumerate}
\item Does the Taylor tower of a functor $F$ converge to anything?
\item Does the Taylor tower converge to $F$?
\end{enumerate}

Information about the first question can be obtained from \refT{homogclass}, the characterization of homogeneous functors, and there is an easy answer if one can compute the connectivity of the derivatives of a functor. The second is much more difficult. This section will first discuss some generalities regarding convergence, including the useful notion of $\rho$-analyticity, where the integer $\rho$ is analogous to a radius of convergence. \refS{convemb} will tackle the convergence question for spaces of embeddings, so the reader has a sense of what types of arguments go into proving convergence results in a specific example.

\subsection{Convergence of the series}

For a smooth function $f:\R\to \R$ with Taylor series $\sum a_k\frac{x^k}{k!}$, the radius of convergence $r$ is the largest value of $r$ such that $\sum a_k\frac{x^k}{k!}$ converges absolutely for $|x|<r$. Thus there are two possibilities for the convergence of the series: either it converges only at $0$, or it converges on an open interval centered at $0$.


We would not speak of convergence of the Taylor series of a functor $F$ unless the homotopy type of $T_kF$ stabilizes with $k$; that is, unless the maps $T_kF\to T_{k-1}F$ have connectivity increasing to infinity with $k$. For a functor $F$, we are interested in the homotopy type of $\holim_k T_kF$, and whether the homotopy type of $T_kF$ ``stabilize'' as $k$ increases. One way to detect this is to study the maps $T_kF\to T_{k-1}F$. If their connectivities increase to infinity with $k$, then we would say that the Taylor series converges, and \refT{homogclass} is useful in giving us a means to attack this. In particular, if the derivatives have increasing connectivity, this will ensure these maps are highly connected. 

\begin{prop}\label{P:seriesconv}
For a good functor $F$, if $F^{(k)}(\emptyset)$ is $c_k$-connected, then $L_kF(M)$ is $(c_k-km)$-connected. More generally, if $U$ has handle dimension $j$, then $L_kF(U)$ is $(c_k-kj)$-connected.
\end{prop}

The homogeneous classification theorem tells us that $L_kF(M)=\hofiber(T_kF(M)\to T_{k-1}F(M))$ is equivalent to the space of compactly supported sections of a fibration over $\binom{M}{k}$ whose fibers are the derivatives $F^{(k)}(\emptyset)$. Thinking of a section space as a twisted mapping space, standard obstruction theory arguments (see \refP{mapconnectivity}) show that if $c_k$ is the connectivity of $F^{(k)}(\emptyset)$, then $L_kF(M)$ is $(c_k-km)$-connected. \refP{mapconnectivity} contains the basic idea. In any case, the Taylor tower of $F$ converges for all $U$ of handle dimension $\leq j$ if $c_k-kj$ tends to infinity with $k$.

We can see that the analog of the radius of convergence has something to do with handle dimension, although we have not yet tackled this in a serious way. This is discussed below as the notion of $\rho$-analyticiy of a functor.

\subsection{Convergence to the functor}

We would certainly say that $T_kF$ converges to $F$ if the canonical map $F\to \holim_kT_kF$ is an equivalence. In this case, the connectivity of $L_kF$ informs us about the connectivity of the ``remainder'' $R_kF=\hofiber(F\to T_kF)$.

\begin{prop}\label{P:fibersconv}
For a good functor $F$, if $F\to\holim_kF$ is an equivalence and $L_{k+1}F$ is $c_k$-connected, where $c_k$ is an increasing function of $k$, then $F\to T_kF$ is $c_k$-connected.
\end{prop}

\begin{proof}
Since $L_{k+1}F=\hofiber(T_{k+1}F\to T_kF)$ is $c_k$-connected, $T_{k+1}\to T_kF$ is $(c_k+1)$-connected, and since $c_k$ is an increasing function of $k$, it follows that $T_lF\to T_kF$ is $(c_k+1)$-connected for all $l>k$. Since $F\to\holim_kT_kF$ is an equivalence, $F\to T_{k}F$ is $c_k$-connected as well.
\end{proof}

Although it may be difficult to establish a homotopy equivalence $F\to \holim_k T_kF$, in practice it is feasible understand the connectivity of $L_kF$ by \refP{seriesconv}, since it reduces to computing the connectivity of the derivatives $F^{(k)}(\emptyset)$. Hence even with a lack of knowledge of convergence, we can formulate conjectures about the connectivities of the maps $F\to T_kF$ based on the connectivity of $L_kF$. Understanding the difference between $F$ and $T_kF$ is a natural question in ordinary calculus as well. We are often interested in the error $R_k(x)=|f(x)-T_kf(x)|$ for certain $x$. For $f$ smooth on $[-r,r]$ and satisfying $|f^{(k+1)}|\leq M_k$ on $(-r,r)$, we have $R_k(x)\leq M_k\frac{r^{k+1}}{(k+1)!}$. If $M_k\frac{r^{k+1}}{(k+1)!}\to 0$ as $k\to\infty$, then we would say that $f$ is analytic on $(-r,r)$; that is, its Taylor series converges to it. We wish, therefore, to answer the following questions:

\begin{enumerate}
\item What is the analog of a radius of convergence?
\item What should we mean by a bound on $f^{(k+1)}$ within the radius of convergence?
\item How can we estimate the ``error'' $R_kF=\hofiber(F\to T_kF)$?
\end{enumerate}

Briefly, the answer to the first question is that the radius of convergence is a positive integer $\rho$. An open set $V$ which is the interior of a smooth compact codimension $0$ submanifold $L$ of $M$ is within the radius of convergence if the handle dimension of $L$ is less than $\rho$. The answer to the second lies in our criticism given in the last section of \refS{derivatives} of our definition of the derivatives of $F$. Our definition of derivative only allows the attaching of a handle of dimension $0$ (disjoint union), while we will need to understand what happens for more general unions. A similar comment applies to higher derivatives. We will expand on all of this below.

To answer the third question, note that we are asking about the extent to which a given functor $F$ fails to be polynomial of degree $\leq k$. We have two options available to us. The first is to study the homotopy fiber of $F(V)\to T_kF(V)$. This has the advantage that it is a natural transformation of functors, and it is the connectivity of this map we are ultimately interested in. Unfortunately, the target is a homotopy limit over a category not very accessible to computation. The other option is to study the extent to which the functor $F$ fails to satisfy the definition of polynomial. This is much more computationally feasible, because it involves values of the original functor on certain kinds of cubical diagrams. 

Suppose $F:\mathcal{O}(M)\to\Top$ is a functor and $\rho>0$ is an integer. For $k>0$, let $P$ be a smooth compact codimension $0$ submanifold of $M$, and $Q_0,\ldots, Q_k$ be pairwise disjoint compact codimension $0$ submanifolds of $M-\interior(P)$. Suppose further that $Q_i$ has handle dimension $q_i<\rho$. Let $U_S=\interior(P\cup Q_S)$.

\begin{defin}
The functor $F$ is \textsl{$\rho$-analytic with excess $c$} if the $(k+1)$-cube $S\mapsto F(U_S)$ is $(c+\sum_{i=0}^k(\rho-q_i))$-cartesian.
\end{defin}

This is the analog of a bound on $f^{(k)}(x)$ for $x$ close to $0$. In this case, close to zero means small handle dimension, and the $(k+1)$-cube $S\mapsto F(U_S)$ certainly resembles a more general $(k+1)$st derivative-like expression. We will see shortly that $\rho$ gives the radius of convergence of the Taylor tower of $F$. Note that this definition is concerned with something close to the $k$th derivative of $F$ at $P$, although we allow ourselves to study multirelative differences not just involving disjoint open balls, but arbitrary manifolds with bounded handle dimension. It is this definition that gives us our answer to the second question above, as we will see in the next theorem, which is the estimate for the error $R_kF=\hofiber(F\to T_kF)$.

\begin{thm}\label{T:rho-analytic}\cite[Theorem 2.3]{GW:EI2}
If $F$ is $\rho$-analytic with excess $c$, and if $U\in\mathcal{O}(M)$ is the interior of a smooth compact codimension $0$ submanifold of $M$ with handle dimension $q<\rho$, then the map $F(U)\to T_{k}F(U)$ is $(c+(k+1)(\rho-q))$-connected.
\end{thm}

\begin{cor}\label{C:rho-analyticcor}\cite[Corollary 2.4]{GW:EI2}
Suppose $F$ is $\rho$-analytic with excess $c$. Then for all $U\in\mathcal{O}(M)$ which are the interior of a compact codimension $0$ submanifold of handle dimension $<\rho$, the map $F(U)\to \holim_k T_kF(U)$  is an equivalence.
\end{cor}

This follows since the connectivities of the maps $F(U)\to T_kF(U)$ increase to infinity with $k$ if the handle dimension of $U$ is less than $\rho$. Thus we see how the handle dimension can be thought of as the radius of convergence, where an open set is measured by its handle dimension.

We will not give the proof of \refT{rho-analytic}, although we would like to make a few remarks. The strategy of the proof is similar to the inductive proof of \refT{charlinear}. We are interested in the connectivity of the map $F(U)\to T_kF(U)$, and as usual, it suffices to study the special case where $U$ is the interior of a smooth compact codimension $0$ submanifold $L$ of $M$. Using a handle decomposition, we select pairwise disjoint closed subsets $A_0,\ldots, A_k$ such that for $S\neq\emptyset$, $U_S=U-\cup_{i\in S}A_i$ is the interior of a compact smooth codimension $0$ submanifold whose handle dimension is strictly less than the handle dimension of $L$. We then consider the diagram
$$
\xymatrix{
F(U)\ar[r]\ar[d]& T_kF(U)\ar[d]\\
\holim_{S\neq \emptyset}F(U_S)\ar[r] & \holim_{S\neq \emptyset}T_kF(U_S)
}.
$$
The right vertical arrow is an equivalence since $T_kF$ is polynomial of degree $\leq k$, and by induction we can get a connectivity estimate for the bottom horizontal arrow. Together these give an estimate for the connectivity of $F(U)\to T_kF(U)$. We have a connectivity for the left vertical arrow by assuming $F$ is $\rho$-analytic. The next section is devoted to understanding how to obtain such connectivity estimates in the case $k=1,2$ for the functor $F(U)=\Emb(U,N)$. In particular, the difficult task is verifying that a given functor is $\rho$-analytic for some $\rho$, which gives a connectivity estimate for the left vertical arrow. Before we embark on this, let us state one more corollary regarding convergence. The next result is the analog of the uniqueness of analytic continuation.

\begin{cor}\label{C:rho-analyticcor2}\cite[Corollary 2.6]{GW:EI2}
Suppose $F_1\to F_2$ is a natural transformation of $\rho$-analytic functors, and that $F_1(U)\to F_2(U)$ is an equivalence whenever $U\in\mathcal{O}_k(M)$ for some $k$. Then $F_1(V)\to F_2(V)$ is an equivalence for all $V$ which are the interior of a smooth compact codimension $0$ submanifold of handle dimension less than $\rho$.
\end{cor}

\begin{proof}
Suppose $V\in\mathcal{O}(M)$. Consider the following diagram.
$$
\xymatrix{
F_1(V)\ar[r]\ar[d]& F_2(V)\ar[d]\\
\holim_kT_kF_1(V)\ar[r] & \holim_kT_kF_2(V)
}.
$$
Since $F_1(U)\to F_2(U)$ is an equivalence whenever $U$ is in $\mathcal{O}_k(M)$ for any $k$, it follows from \refT{polychar} that $T_kF_1\to T_kF_2$ is an equivalence for all $k$. Hence the lower horizontal arrow is an equivalence for all $V$. If the handle dimension of $V$ is less than $\rho$, then $F_1(V)\to\holim_kT_kF_1(V)$ and $F_2(V)\to\holim_kT_kF_2(V)$ are equivalences by \refC{rho-analyticcor}, so $F_1(V)\to F_2(V)$ is an equivalence.
\end{proof}

\section{Convergence for Spaces of Embeddings}\label{S:convemb}

The following is a theorem due to Klein and Goodwillie about the convergence of the Taylor tower of the embedding functor. A version for spaces of Poincar\'e embeddings has appeared in \cite{GK:PE}, which is an important step in proving the result below, which will appear in \cite{GK:SE}.

\begin{thm}\label{T:GK}
The functor $U\mapsto\Emb(U,N)$ is $n-2$ analytic with excess $3-n$. Hence, if $M$ is a smooth closed manifold of dimension $m$, and $N$ a smooth manifold of dimension $n$, then the map
$$
\Emb(M,N)\longrightarrow T_k\Emb(M,N)
$$
is $[k(n-m-2)+1-m]$-connected. In particular, if $n-m-2>0$, then the canonical map
$$
\Emb(M,N)\longrightarrow\holim_{k}T_k\Emb(M,N)
$$
is a homotopy equivalence.
\end{thm}

The proof of this theorem goes beyond the scope of this work, although we wish to present  some of the ideas involved in arriving at such estimates. Note that the estimate for the map $\Emb(M,N)\to T_k\Emb(M,N)$ can be conjectured using \refP{fibersconv}; we will compute the connectivity of the derivatives of embeddings below. Note also that in the case $m=1$ and $n=3$ (essentially knot theory), we do not have convergence (although the theorem still gives a non-trivial answer).

One can obtain the connectivity of $\Emb(M,N)\to T_1\Emb(M,N)$ ``by hand'' without too much work, and some of the ideas that go into one version of this computation (the second proof of \refT{etotone} below) are important in obtaining estimates for all $k$. We will also discuss a weaker estimate for the map $\Emb(M,N)\to T_2\Emb(M,N)$. The techniques required for the results above are far beyond the scope of this work, and involves important relationships between embeddings, pseudoisotopies, and diffeomorphisms, as well as some surgery theory.

\subsection{Connectivity of the derivatives of embeddings}

The first step in understanding some of the ideas that go into establishing the analyticity of the embedding functor is to compute the connectivity of the derivatives of the embedding functor.

\begin{thm}
Let $U=\coprod_i B_i\subset M$ be a disjoint union of $k$ open balls. For $S\subset\underline{k}$, let $U_S=U-\cup_{i\in S} B_i$. The $k$-cube $S\mapsto\Emb(U_S,N)$ is $((k-1)(n-2)+1)$-cartesian. That is if $E(U)=\Emb(U,N)$, then $E^{(k-1)}(\emptyset)$ is $(k-1)(n-2)$-connected.
\end{thm}

Let us begin with an observation that will simplify things. For a subset $S$ of $\underline{k}$, the projection map $\prod_{i\notin S}B_i\times\Emb(U_S,N)\to\Emb(U_S,N)$ is an equivalence because balls and products of balls are contractible (if $S=\underline{k}$, we take $\prod_{i\notin S}B_i$ to be a point). Let $C(j,N)$ denote the configuration space of $j$ points in $N$. The map $\prod_{i\notin S}B_i\times\Emb(U_S,N)\to C(k-|S|,N)\times\Imm(U_S,N)$ which is induced by the map which sends $((x_1,\ldots, x_k),f)$ to $((f(x_1),\ldots, f(x_k),(df_{x_1},\ldots, df_{x_k}))$ is an equivalence for all $S$ (where again the product of balls is taken to be a point if $S=\underline{k}$). Hence $S\mapsto\Emb(U_S,N)$ is $K$-cartesian if and only if $S\mapsto C(k-|S|,N)\times \Imm(U_S,N)$ is $K$-cartesian. The cube $S\mapsto \Imm(U_S,N)$ is homotopy cartesian whenever $k\geq 2$ because $\Imm(-,N)$ is polynomial of degree $\leq 1$. Therefore $S\mapsto\Emb(U_S,N)$ is $K$-cartesian if and only if $S\mapsto C(k-|S|,N)$ is $K$-cartesian for $k\geq 2$. For illustration, we will only prove this in the case where $k=2$. The cases $k\geq 3$ are straightforward enough, and all that they require is an application of the Blakers-Massey \refT{BM}.

\begin{proof}
For $k=2$, we are looking at the square
$$
\xymatrix{
\Emb(U,N)\ar[r]\ar[d]& \Emb(U_0,N)\ar[d]\\
\Emb(U_1,N)\ar[r] & \Emb(\emptyset,N)
}.
$$
By the remarks preceding the proof, this square is $K$-cartesian if and only if 
$$
\xymatrix{
C(2,N)\ar[r]\ar[d]& C(1,N)\ar[d]\\
C(1,N)\ar[r] & C(\emptyset,N)
}
$$
is $K$-cartesian. The maps in this diagram are fibrations, and taking fibers vertically over $p\in C(1,N)$ yields the $1$-cube $N-\{p\}\to N$, which is an $(n-1)$-connected map, and hence the original square is $(n-1)$-cartesian by \refP{cartesianfibercube}.
\end{proof}

As we mentioned, the Blakers-Massey \refT{BM} needs to be applied for higher $k$. For instance, the case $k=3$ ends with fibering over $(p,q)\in C(2,N)$ and observing that the square
$$
\xymatrix{
N-\{p,q\}\ar[r]\ar[d]& N-p\ar[d]\\
N-q\ar[r] & N
}
$$
is a homotopy pushout and is $(2n-3)$-cartesian by the Blakers-Massey \refT{BM}.

\subsection{Connectivity estimates for the linear and quadratic stages for embeddings}

We will give two proofs of the following theorem. The second requires a disjunction result from the next section, but beyond this, it is almost identical to the proof of \refT{charlinear}.

\begin{thm}\label{T:etotone}
The map $\Emb(M,N)\to T_1\Emb(M,N)$ is $(n-2m-1)$-connected. In fact, if $V\subset M$ is the interior of a compact codimension $0$ handlebody with handle dimension $k$, then the map is $(n-2k-1)$-connected.
\end{thm}

The first proof is much easier and employs general position arguments, although it only gives the connectivity estimate in terms of the dimension of $M$, not the stronger statement involving the handle dimension. The second uses a bit more machinery, but reduces the proof to the special case where $M$ is the disjoint union of balls via an induction argument on the handle dimension, but requires a disjunction result from the next section. Hopefully this further convinces the reader of the importance of derivatives. Its methods are also important in organizing the proof of the connectivity estimate for $\Emb(M,N)\to T_j\Emb(M,N)$ for all $j$.

\begin{proof}[First Proof]
We have already mentioned that $T_1\Emb(M,N)\simeq\Imm(M,N)$. Let $h:S^k\to\Imm(M,N)$ be a map with adjoint $H:M\times S^k\to N$. Consider the map $\tilde{H}:M\times M\times S^k \to N\times N$ defined by $\tilde{H}(x,y,s)=(H(x,s),H(y,s))$. We can arrange, by a small homotopy, for $H$ to be smooth and $\tilde{H}$ to be transverse to the diagonal. Let $D=\tilde{H}^{-1}(\Delta(N))$ be the inverse image of the diagonal. It is a submanifold of $M\times M\times S^k$ of dimension $2m+k-n$, which is empty if $k<n-2m$, and in this case, the map $h$ clearly has image in $\Emb(M,N)$. A similar argument shows that a homotopy $h:S^k\times I\to\Imm(M,N)$ lifts to $\Emb(M,N)$ if $k<n-2m-1$, and it follows that the inclusion $\Emb(M,N)\to\Imm(M,N)$ is $(n-2m-1)$-connected.
\end{proof}

\begin{proof}[Second Proof]
We will induct on $k$. For the base case $k=0$, let $l$ be the number of components of $V$. The result is trivial, and the map in question is an equivalence, when $l=0,1$. Suppose $l\geq 2$. Consider the sequence
$$
\Emb(V,N)\to T_l\Emb(V,N)\to T_{l-1}\Emb(V,N)\to\cdots\to T_1\Emb(V,N).
$$
The map $\Emb(V,N)\to T_l\Emb(V,N)$ is an equivalence since $V$ is a final object in $\mathcal{O}_l(V)$. By the classification \refT{homogclass} of homogeneous functors, we have that $L_j\Emb(V,N)=\hofiber(T_j\Emb(V,N)\to T_{j-1}\Emb(V,N))$ is equivalent to $\Gamma^c(\binom{V}{j},\Emb^{(j)}(\emptyset))$. Since $V$ has handle dimension $0$, $\binom{V}{j}$ also has handle dimension $0$, and the fibers are thus $(j-1)(n-2)$-connected; in other words, the map $T_j\Emb(V,N)\to T_{j-1}\Emb(V,N)$ is $((j-1)(n-2)+1)$-connected. This is true no matter what basepoint is chosen, provided $m<n$. It follows that the composed map $\Emb(V,N)\to T_1\Emb(V,N)$ is $(n-1)$-connected.

Now suppose $k>0$. Let $V=\interior(L)$. For $j=1$ to $s$, let $e_j: D^{k}\times D^{n-k}\to L$ denote each of the $s$ $k$-handles. Assume $e_j^{-1}(\del L)=\del D^k\times D^{n-k}$ for all $j$. Since $k>0$, we may choose pairwise disjoint closed disks $D_0,D_1$ in the interior of $D^k$, and put $A^j_i=e_j(D_i\times D^{n-k})\cap V$. Then each $A_i^j$ is closed in $V$, and if we set $A_i=\cup_{j=1}^s A_i^j$, then for each nonempty subset $S$ of $\{0,1\}$, $V_S=V-\cup_{i\in S}A_i$ is the interior of a smooth compact codimension $0$ submanifold of $M$ of handle dimension strictly less than $k$. See Figure \ref{handles} for a low-dimensional picture where there are four disks $D_i$ instead of just two.

\begin{figure}
\begin{center}
\input{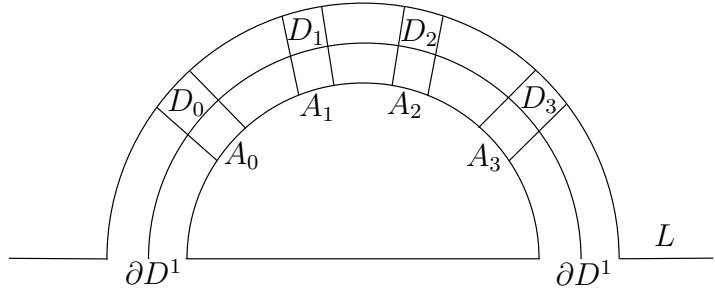}
\caption{A picture of four disks $D_i$ in the core of a $1$-handle $D^1\times D^1$ attached to $L$ along $\del D^1\times D^1$ and their corresponding thickenings $A_i$. The $D_i$ are subsets of the core $D^1\times\{0\}$ (which itself is depicted as the curve in the middle of the handle), and $A_i=D_i\times D^1\subset D^1\times D^1$. Note that removing $k\geq 1$ of the $A_i$ leaves a manifold with $(k-1)$ extra $0$-handles, but one fewer $1$-handle}
\label{handles}
\end{center}
\end{figure}

In the following square diagram,
$$
\xymatrix{
\Emb(V,N)\ar[r]\ar[d]& T_1\Emb(V,N)\ar[d]\\
\holim_{S\neq \emptyset}\Emb(V_S,N)\ar[r] & \holim_{S\neq \emptyset}T_1\Emb(V_S,N)
}
$$
the right vertical arrow is again an equivalence because $T_1\Emb(-,N)$ is polynomial of degree $\leq 1$, and by induction for all $S\neq\emptyset$, $\Emb(V_S,N)\to T_1\Emb(V_S,N)$ is $(n-2(k-1)-1)$-connected, and by \refP{ScubeTcartesian}, the map of homotopy limits has connectivity $n-2(k-1)-1-2+1=n-2k$. By \refT{tonedisj}, the left vertical map is $(n-2k-1)$-connected, and it follows that the top horizontal map is $(n-2k-1)$-connected.
\end{proof}

The base case of the induction on handle dimension above required an argument which was different than the inductive step. In particular it required knowledge of the higher derivatives, and we do not see a way around this. Attempts to mimic the inductive step for the base case yields connectivity estimates which are less than those desired.

\begin{thm}\label{T:etottwo}
The map $\Emb(M,N)\to T_2\Emb(M,N)$ is $(2n-3m-3)$-connected. In fact, if $V\subset M$ is the interior of a compact codimension $0$ submanifold of $M$ whose handle dimension is $k$, then the map $\Emb(V,N)\to T_2\Emb(V,N)$ is $(2n-3k-3)$-connected.
\end{thm}

The second proof of \refT{etotone} can be adapted with very few changes. The only changes (besides the connectivity estimates themselves) are that the pairwise disjoint closed subsets chosen are three in number, and instead of referencing \refT{tonedisj}, we reference \refT{ttwodisj}. However, \refT{ttwodisj} is weaker than what we need, and we can really only claim to prove a weaker version of \refT{etottwo}. The issue here is that there is a stronger version of \refT{ttwodisj} which we are unable to prove by elementary means.

\begin{thm}\label{T:etottwoweak}
With hypotheses as in \refT{etottwo}, the map $\Emb(V,N)\to T_2\Emb(V,N)$ is $(2n-4k-3)$-connected.
\end{thm}

\subsection{Some disjunction results for embeddings}\label{S:disj}

For the second proof of \refT{etotone} we needed an estimate for how cartesian the square $\mathcal{E}$
$$
\xymatrix{
\Emb(V,N)\ar[r]\ar[d]& \Emb(V_0,N)\ar[d]\\
\Emb(V_1,N)\ar[r] & \Emb(V_{01},N)
}
$$
is. Here $V=V_\emptyset$ is the interior of some smooth compact codimension $0$ submanifold of $M$ with handle dimension $k$, and, for $S\neq\emptyset$, the $V_S$ are the interiors of compact codimension $0$ submanifolds of handle dimension less than $k$. As in the proof of \refT{etotone}, let $V=\interior(L)$. We chose each $A_i$ to be a (union of) the product of a $k$-dimensional disk with an $(m-k)$-dimensional disk. Note that $L_S=L-\cup_{i\in S}A_i$ is \emph{not} compact, but its interior is the interior of a smooth compact codimension $0$ submanifold of $M$. This is important to note because below we will work not with the open sets that appear in $\mathcal{E}$, but with their closed counterparts $L$ and the $A_i$.

Let us first consider a formally similar situation. Suppose $Q_0$ and $Q_1$ are smooth closed manifolds of dimensions $q_0$ and $q_1$ respectively, and let $Q_S=\cup_{i\notin S}Q_i$ for $S\subset\{0,1\}$. Consider the square $\mathcal{S}=S\mapsto\Emb(Q_S,N)$:
$$
\xymatrix{
\Emb(Q_0\cup Q_1,N)\ar[r]\ar[d]& \Emb(Q_0,N)\ar[d]\\
\Emb(Q_1,N)\ar[r] & \Emb(\emptyset,N).
}
$$
It is enough by \refP{cartesianfibercube} to choose a basepoint in $\Emb(Q_0\cup Q_1,N)$ and take fibers vertically and compute the connectivity of the map of homotopy fibers. By the isotopy extension theorem, the map $\Emb(Q_0\cup Q_1,N)\to\Emb(Q_1,N)$ is a fibration with fiber $\Emb(Q_0, N-Q_1)$. In \refT{tonedisj} we will show that $\Emb(Q_0,N-Q_1)\to\Emb(Q_0,N)$, and hence the square $\mathcal{S}$, is $(n-q_0-q_1-1)$-cartesian. Although the squares $\mathcal{E}$ and $\mathcal{S}$ are formally similar, it is not clear how to use \refT{tonedisj} to give an estimate for how cartesian the square $\mathcal{E}$ is. 

First note that we can generalize the situation in the square $\mathcal{S}$ to a relative setting. That is, suppose $Q_0$, $Q_1$ and $N$ have boundary, and that embeddings $e_i:\del Q_i\to\del N$ have been selected to have disjoint images. Let $\Emb_\del(Q_S,N)$ be the space of embeddings $f:Q_S\to N$ such that the restriction of $f$ to $\del Q_S$ is equal to $e_S$, and such that $f^{-1}(\del N)=\del Q_S$. Then it is also true that 
$$
\xymatrix{
\Emb_\del(Q_0\cup Q_1,N)\ar[r]\ar[d]& \Emb_\del(Q_0,N)\ar[d]\\
\Emb_\del(Q_1,N)\ar[r] & \Emb_\del(\emptyset,N)
}
$$
is $(n-q_0-q_1-1)$-cartesian; in particular, the proof of this is identical to that of \refT{tonedisj} with the exception of repeating the phrase ``relative to the boundary'' over and over. 

We can make a further generalization to the case of compact manifold triads (defined in \refS{convs}). Suppose the $Q_i$ are compact $n$-dimensional manifold triads of handle dimension $q_i$, where $n-q_i\geq 3$, and $Y$ is an $n$-dimensional smooth manifold with boundary. In this case embeddings $e_i:\del_0Q_i\to\del N$ have been chosen, and we let $\Emb_{\del_0}(Q_S,N)$ stand for the obvious thing.

\begin{thm}\label{T:disj1}\cite[Theorem 1.1]{GW:EI2}
$$
\xymatrix{
\Emb_{\del_0}(Q_0\cup Q_1,N)\ar[r]\ar[d]& \Emb_{\del_0}(Q_0,N)\ar[d]\\
\Emb_{\del_0}(Q_1,N)\ar[r] & \Emb_{\del_0}(\emptyset,N)
}
$$
is $(n-q_0-q_1-1)$-cartesian.
\end{thm}

This can be generalized to the case where the dimension of the $Q_i$ is $m\leq n$, essentially by a thickening of the $m$-dimensional $Q_i$ by the disk bundle of an $(n-m)$-plane bundle. 

\begin{prop}\label{P:disj1dim}\cite[Observation 1.3]{GW:EI2}
\refT{disj1} is true if $\dim(Q_i)=m\leq n$.
\end{prop}

The rough idea of the proof is to assume that $Y$ is embedded in $\R^{n+k}$ and let $Gr_{n-m}=\colim_{k}Gr_{n-m+k}(\R^{n+k})$ be a limit of Grassmannians. Consider the map $\Emb(Q_S,Y)\to\Map(Q_S,Gr_{n-m})$ given by assigning an embedding $f$ to its normal bundle $\nu_f$. The homotopy fiber of this map over some $\eta$ can be identified with the space of embeddings of the disk bundle of $\eta$ over $Q_S$. Since $S\mapsto\Map(Q_S,Gr_{n-m})$ is homotopy cartesian (because $\Map(-,X)$ is polynomial of degree $\leq 1$), by \refP{cartesianfibercube}, the square of homotopy fibers is $(n-q_0-q_1-1)$-cartesian if and only if the square $S\mapsto\Emb_\del(Q_S,Y)$ is $(n-q_0-q_1-1)$-cartesian. Note, however, that this introduces more corners, since the closed disk bundle of a smooth manifold with boundary is already a compact manifold triad itself. The new corners due to the disk bundle are introduced along the corner set of the original compact manifold triad. It will do no harm to ignore this.

Without any changes whatsoever we can assume the $Q_i$ are submanifolds of an $m$-dimensional manifold $M$. Now we are in a position to describe a situation which is directly related to the square $\mathcal{E}$, and we generalize this situation further by introducing a new manifold $P$. Suppose that $P$ is a smooth compact codimension $0$ manifold triad in $M$, $Q_0,Q_1$ are smooth compact codimension $0$ manifold triads in $M-\interior(P)$, and that the handle dimension of $Q_i$ satisfies $n-q_i\geq 3$. Put $Q_S=\cup_{i\notin S}Q_i$.

\begin{prop}\label{P:embcartesianfiber}
The square $S\mapsto\Emb(P\cup Q_S,N)$ is $(n-q_0-q_1-1)$-cartesian.
\end{prop}

\begin{proof}
The square $S\mapsto\Emb(P,N)$ is homotopy cartesian since all maps are equivalences, and hence $S\mapsto\Emb(P\cup Q_S,N)$ is $(n-q_0-q_1-1)$-cartesian if and only if the square of homotopy fibers $S\mapsto\hofiber(\Emb(P\cup Q_S,N)\to\Emb(P,N))$ is $(n-q_0-q_1-1)$-cartesian. The map $\Emb(P\cup Q_S,N)\to\Emb(P,N)$ is a fibration with fiber $\Emb_{\del_0}(Q_S,N-P)$, which is $(n-q_0-q_1-1)$-cartesian by \refT{disj1}.
\end{proof}

We finally arrive at the technical statement which relates the open sets in square $\mathcal{E}$ with the closed sets we have been considering.

\begin{cor}\label{disj1open}\cite[Corollary 1.4]{GW:EI2}
Let $P, Q_0,Q_1$ be as in \refP{embcartesianfiber}, and set $V_S=\interior(P\cup Q_S)$. Then $S\mapsto\Emb(V_S,N)$ is $(n-q_0-q_1-1)$-cartesian.
\end{cor}

To connect this explicitly with the square $\mathcal{E}$, we choose the $Q_i$ to be the $A_i$ considered in \refT{etotone}, and $P$ to be the closure of $L-(A_0\cup A_1)$. We now proceed to give the promised disjunction results.

\begin{thm}\label{T:tonedisj}
Suppose $P$ and $Q$ are smooth compact submanifolds of an $n$-dimensional manifold $N$ of dimensions $p$ and $q$ respectively. The inclusion map $\Emb(P,N-Q)\to\Emb(P,N)$ is $(n-p-q-1)$-connected.
\end{thm}

An important special case is when both $P$ and $Q=\ast$ are points, which says that $N-\ast\to N$ is $(n-1)$-connected. The rough idea, expanded in the proof below, is that any map $S^k\to N$ misses a point if $k<n$, and that the same is true of any homotopy $S^k\times I$ if $k<n-1$. The former proves the map of homotopy groups is surjective if $k<n$ and the latter that it is injective if $k<n-1$.

\begin{proof}
We will not fuss about basepoints. The following argument can be adapted to accomodate them. Let $S^k\to\Emb(P,N)$. We may regard this as a map $S^k\times P\to N$, and by a small homotopy we can make it both smooth and transverse to $Q\subset N$. If $k+p<n-q$, equivalently, $k<n-p-q$, transverse intersection means empty intersection, and hence we have a map $S^k\times P\to N-Q$, which gives us our desired map $S^k\to\Emb(P,N-Q)$. A similar argument shows that any homotopy $S^k\times I\to\Emb(P,N)$ lifts to $\Emb(P,N-Q)$ if $k<n-p-q-1$, hence the map in question is $(n-p-q-1)$-connected.
\end{proof}


We can piggyback on the previous result to obtain the following generalization.

\begin{thm}\label{T:ttwodisj}
The square of inclusion maps
$$\xymatrix{
\Emb(P,N-(Q_0\cup Q_1))\ar[r]\ar[d]& \Emb(P,N-Q_0)\ar[d]\\
\Emb(P,N-Q_1)\ar[r] &\Emb(P,N)
}
$$
is $(2n-2p-q_1-q_2-3)$-cartesian.
\end{thm}

Once again the special case where $P$, $Q_0=\ast_0$ and $Q_1=\ast_1$ are all points is a good one to consider before embarking on the proof. In that case, we claim that the square of inclusion maps
$$\xymatrix{
N-(\ast_0\cup \ast_1)\ar[r]\ar[d]& N-\ast_0\ar[d]\\
N-\ast_1\ar[r] & N
}
$$
is $(2n-3)$-cartesian. The square is clearly a homotopy pushout, and since the maps $N-(\ast_0\cup \ast_1)\to N-\ast_i$ are $(n-1)$-connected for $i=0,1$, by the Blakers-Massey Theorem, the square is $(2n-3)$-cartesian. Unfortunately, in the general case the square will not be a homotopy pushout, but it is close to being one, and we will use a generalization of the Blakers-Massey Theorem to complete the proof.

\begin{proof}
We claim that the square in question is $(2n-2p-q_0-q_1-1)$-cocartesian (see \refD{totalcofiber} in the Appendix). Given this, since the maps $\Emb(P,N-(Q_0\cup Q_1))\to\Emb(P,N-Q_i)$ are $(n-p-q_j-1)$-connected for $i=0,1$ and $j=1,0$ respectively by \refT{tonedisj}, it follows by the generalized Blakers-Massey \refT{BMgen} that the square is $\min\{2n-2p-q_0-q_1-1,2n-2p-q_0-q_1-3\}$-cartesian, and the result follows.

Now we will establish the estimate for how cocartesian this square is. Let $h:S^k\to\Emb(P,N)$ be a map, and $H:P\times S^k\to N$ its adjoint. We wish to construct a lift to $\hocolim(\Emb(P,N-Q_0)\leftarrow \Emb(P,N-(Q_0\cup Q_1))\rightarrow \Emb(P,N-Q_1))$.  It is enough if $H$ has the property that for all $s\in S^k$, $H(P\times\{s\})\subset N-Q_i$ for some $i$ (possibly both).

By a small homotopy, we can make the map $H:P\times S^k\to N$ smooth and transverse to $Q_0$ and $Q_1$. For $i=0,1$, let $W_i=H^{-1}(Q_i)$. Then $W_0$ and $W_1$ are disjoint submanifolds of $P\times S^k$ of dimensions $p+k+q_0-n$ and $p+k+q_1-n$ respectively. Consider the map $d:W_0\times W_1\to S^k\times S^k$ induced by projection to $S^k$. Again by a small homotopy of $H$ we can arrange for this map to be transverse to the diagonal, and let $D=d^{-1}(\Delta(S^k))$. The dimension of $D$ is $2p+q_0+q_1-2n+k$, which is negative (and hence $D$ empty) if $k<2n-2p-q_0-q_1$, so that $h$ lifts if $k\leq 2n-2p-q_0-q_1-1$. A similar argument shows a homotopy $h:S^k\times I\to\Emb(P,N)$ lifts if $k\leq 2n-2p-q_0-q_1-2$, and hence the square is $(2n-2p-q_0-q_1-1)$-cocartesian.
\end{proof}


\begin{rem}
\refT{ttwodisj} tells us that the $3$-cube $S\mapsto(\Emb(P\cup Q_S,N)\to\Emb(Q_S,N))$ is $(2n-2p-q_0-q_1-3)$-cartesian. Relabel and put $Q_2=P$ and $q_2=p$, and for $R\subset \{0,1,2\}$ write $R\mapsto\Emb(Q_R,N)$, where $Q_R=\cup_{i\notin R}Q_i$. The estimate $2n-2q_2-q_0-q_1-3$ is not symmetric in the $q_i$; that is, it depended upon us choosing a way to view this $3$-cube as a map of squares. We could therefore slightly improve our result by letting $P=Q_i$, where $Q_i$ minimizes the handle dimension among $Q_0,Q_1,Q_2$. Although this is a slight improvement, it is not the best possible, which is probably not surprising in light of the asymmetry in the $q_i$. One can prove that the $3$-cube $R\mapsto\Emb(Q_R,N)$ is $(2n-q_0-q_1-q_2-3)$-cartesian.
\end{rem}

Let us end by remarking that a lack of convergence does not necessarily mean the Taylor tower does not contain anything interesting. In fact, for embeddings of $I$ in $\R^2\times I$ relative to the boundary (where \refT{GK} doesn't give convergence because the codimension is $2$), \cite{V:FTK} shows that the Taylor series for the embedding functor contains finite type invariants of knots. On a related note, one can study multivariable functors such as $(U,V)\mapsto\Link(U,V;N)$. Here $U$ and $V$ are open subsets of smooth closed manifolds $P$ and $Q$, and $\Link(U,V;N)$ is the space of ``link maps'' $U\to N$, $V\to N$ whose images are disjoint. It is not known whether its (multivariable) Taylor series converges to it, but it is clear that its polynomial approximations are interesting, since, for example, $\hofiber(\Link(P,Q;N)\to T_1\Link(P,Q;N))$ contains the information necessary to define the generalized linking number. See \cite{M:LinkNumber} and \cite{GM:LinksEstimates}.


\section{Appendix}\label{S:app}

This section contains a collection of facts about homotopy limits and cubical diagrams. It is by no means exhaustive, even for the purposes of this paper. What we have included is most of what one needs to be able to understand this work, and we have tried to include only those facts which seem more widely useful in calculus of functors. Proofs are generally omitted, and if given, are very sketchy and are just meant to outline the major ideas and/or given an intuitive understanding. As many references as possible are given.

\subsection{Homotopy limits and colimits}\label{S:holims}

The standard reference for homotopy limits and colimits is \cite{BK}. Others include \cite{Dwy:HT&MC}, \cite{Dwy:HT&CS}, \cite{Hir:MC&loc}, and \cite{Vogt:CommutingHolims}. We begin with an explicit description of the homotopy limit and homotopy colimit. We assume the reader is familiar with simplicial sets and their realizations, as well as over/under categories. We present these models because they involve expressions which are easy for a beginner to grasp, and while other models for homotopy limits and colimits are better for other purposes, these are good for getting  one's hands dirty with categories with finitely many objects and morphisms. For a category $\mathcal{D}$, we write $|\mathcal{D}|$ for the realization of its nerve.

\begin{defin}\label{D:holim}\cite[XI.3.2 and XII.2.1]{BK}
Let $\mathcal{C}$ be a small category and $F:\mathcal{C}\to\Top$ a covariant functor. The \textsl{homotopy limit of $F$}, denoted $\holim_\mathcal{C}F$, is 
$$
\lim(\underset{c}{\prod}\Map(|\mathcal{C}\downarrow c|, F(c))\rightarrow\underset{c\to c'}{\prod}\Map(|\mathcal{C}\downarrow c|, F(c'))\leftarrow\underset{c'}{\prod}\Map(|\mathcal{C}\downarrow c'|, F(c'))).
$$
Dually, the \textsl{homotopy colimit of $F$}, denoted $\hocolim_{C}F$, is
$$
\colim(\underset{c}{\coprod}|c\downarrow\mathcal{C}|\times F(c)\leftarrow\underset{c\to c'}{\coprod}|c'\downarrow\mathcal{C}|\times F(c)\rightarrow\underset{c'}{\coprod}|c'\downarrow\mathcal{C}|\times F(c')).
$$
\end{defin}

The maps in these diagrams are induced by the identity map and functoriality of various functors. That is, $c\to c'$ induces maps $F(c)\to F(c')$, $|c'\downarrow\mathcal{C}|\to|c\downarrow\mathcal{C}|$, and $|\mathcal{C}\downarrow c|\to|\mathcal{C}\downarrow c'|$. Note that $\Map(-,-)$ is contravariant in the first variable and covariant in the second. One nice feature of this definition is that it only requires the reader to understand limits and colimits of very simple diagrams, namely $\lim(X_1\to X_{12}\leftarrow X_2)$ and $\colim(X_1\leftarrow X_\emptyset\to X_2)$, which are the fiber product of $X_1$ with $X_2$ over $X_{12}$ and the union of $X_1$ with $X_2$ along $X_\emptyset$ respectively.

The ordinary limit (indirect limit) and colimit (direct limit) of $F$ can be defined by replacing all realizations of over/under categories above with a point. Indeed, the limit is a subspace of the product $\prod_{c\in\mathcal{C}}F(c)$ and the colimit is a subspace of the coproduct $\coprod_{c\in\mathcal{C}}F(c)$.

An important special case is when $\mathcal{C}=\mathcal{P}(\underline{2})$ is the poset of subsets of $\{1,2\}$. Let $\mathcal{X}:\mathcal{P}(\{1,2\})\to\Top$ be a covariant functor, and write $X_S$ in place of $\mathcal{X}(S)$. Let $\mathcal{P}_0(\underline{2})$ be the subposet of nonempty subsets, and let $\mathcal{P}_1(\underline{2})$ be the subposet of proper subsets. We depict the diagram of spaces as follows.

$$\xymatrix{
X_\emptyset \ar^{f_1}[r]\ar_{f_2}[d]& X_1\ar^{g_1}[d]\\
X_2\ar_{g_2}[r] & X_{12}
}
$$

Although it is somewhat tedious, it is straightforward to show the following from the definition of homotopy limit and homotopy colimit.

\begin{prop}\label{P:holimsquare}
We have
$$
\underset{\mathcal{P}_0(\underline{2})}{\holim}\;F=\{(x_1,\gamma, x_2)\in X_1\times \Map(I,X_{12})\times X_2|\gamma(0)=g_1(x_1), \gamma(1)=g_2(x_2)\}
$$
and
$$
\underset{\mathcal{P}_1(\underline{2})}{\hocolim}\;F=(X_1\coprod X_\emptyset\times I\coprod X_2)/\sim
$$
where for $x\in X_\emptyset$, $(x,0)\sim f_1(x)$ and $(x,1)\sim f_2(x)$. In particular, if $X_2=\ast$ is the one-point space, then $\holim_{\mathcal{P}_0(\underline{2})}F$ is the homotopy fiber of $g_1$ over $g_2(\ast)$, and $\hocolim_{\mathcal{P}_1(\underline{2})}F$ is the homotopy cofiber of $f_1$.
\end{prop}

Thus homotopy (co)fiber is a special case of homotopy (co)limit. It is also true that homotopy (co)limits commute.

\begin{thm}\label{T:holimcommute}
Suppose $F:\mathcal{C}\times\mathcal{D}\to\Top$ is a bifunctor. Then there are homeomorphisms
$$
\underset{\mathcal{C}}{\holim}\;\underset{\mathcal{D}}{\holim}\; F\cong \underset{\mathcal{C}\times\mathcal{D}}{\holim}\; F\cong\underset{\mathcal{D}}{\holim}\;\underset{\mathcal{C}}{\holim}\; F
$$
and
$$
\underset{\mathcal{C}}{\hocolim}\;\underset{\mathcal{D}}{\hocolim}\; F\cong \underset{\mathcal{C}\times\mathcal{D}}{\hocolim}\; F\cong\underset{\mathcal{D}}{\hocolim}\;\underset{\mathcal{C}}{\hocolim}\; F.
$$
\end{thm}

In particular, homotopy limits commute with homotopy fibers and homotopy colimits commute with homotopy cofibers. The following theorem establishes the homotopy invariance of homotopy limits and colimits.

\begin{thm}\label{T:holimheq}\cite[XI.5.6 and XII.4.2]{BK}
Suppose $F\to G$ is a natural transformation of functors from $\mathcal{C}$ to $\Top$. If $F(c)\to G(c)$ is an equivalence for all $c\in\mathcal{C}$, then this induces equivalences
$$
\underset{\mathcal{C}}{\holim}\;F\longrightarrow \underset{\mathcal{C}}{\holim}\;G
$$
and
$$
\underset{\mathcal{C}}{\hocolim}\;F\longrightarrow \underset{\mathcal{C}}{\hocolim}\;G.
$$ 
\end{thm}

This fact is not true for ordinary limits and colimits; one can view the construction of the homotopy limit/colimit as a way to remedy this. Nevertheless, there are conditions under which the categorical (co)limit is equivalent to the homotopy (co)limit; we will not pursue this here, but point out that we have encountered this situation already for square diagrams (where it was enough if the square were a categorical pushout/pullback and a map from the initial object/to the final object were a cofibration/fibration).

\begin{thm}\label{T:holiminitial}\cite[XI.4.1 and XII.3.1]{BK}
Suppose $F:\mathcal{C}\to \Top$ is covariant. If $\mathcal{C}$ has an initial object $c_i$, then $\holim_\mathcal{C}F\simeq F(c_i)$. If $\mathcal{C}$ has a final object $c_f$, then $\hocolim_\mathcal{C}F\simeq F(c_f)$.
\end{thm}

If $F$ is contravariant, then we need to switch ``initial'' with ``final'' in the above statement. The corresponding facts about ordinary limits and colimits are obvious if one defines such notions in terms of universal properties. The following is a useful result we used in the proof of \refT{charlinear}, and is also central to the proof of \refT{homogclass}. It describes a close relationship between homotopy limits and homotopy colimits when the functor $F$ is especially well-behaved. It has a similar flavor to Quillen's Theorems A and B.

\begin{thm}\label{T:holimquasifiber}\cite{Dwy:BG}
If $F:\mathcal{C}\to\Top$ takes all morphisms to homotopy equivalences, then $\hocolim_{\mathcal{C}}F$ quasifibers over $|\mathcal{C}|$, and the space of sections of the associated fibration is equivalent to $\holim_\mathcal{C}F$.
\end{thm}

The quasifibration statement is at least relatively easy to believe. If we let $\ast:\mathcal{C}\to\Top$ denote the functor which takes all objects to the one-point space, then $\hocolim_{\mathcal{C}}\ast\simeq|\mathcal{C}|$, and the natural transformation $F\to\ast$ induces the map $\hocolim_\mathcal{C}F\to|\mathcal{C}|$. Since $F$ takes all morphisms to homotopy equivalences, the fibers $F(c)$ all have the same homotopy type.

\subsection{The functor $\Map(-,-)$}

We can regard $\Map(-,-):\Top\times\Top\to\Top$ as a bifunctor which is contravariant in the first variable and covariant in the second variable. We are mostly interested in $\Map(-,Z)$ and its variants for a fixed $Z$, but it is also useful to consider $\Map(X,-)$ for a fixed $X$.

\begin{prop}\label{P:mapconnectivity}
For a finite complex $K$ of dimension $k$, the functor $\Map(K,-)$ takes $j$-connected maps to $(j-k)$-connected maps. In particular, if $Z$ is $j$-connected, then $\Map(K,Z)$ is $(j-k)$-connected.
\end{prop}

This can be proved by standard obstruction theory arguments. It is true of more general mapping spaces to, such as the space of sections of a fibration $p:E\to B$, which we loosely think about as the space of maps from $B$ into the fiber of $p$.

\begin{prop}\label{P:maphocolimtoholim}\cite[XII.4.1]{BK}
The functor $\Map(-,X):\Top\to\Top$ takes (homotopy) colimits to (homotopy) limits. That is, if $\mathcal{C}$ is a small category and $\mathcal{X}:\mathcal{C}\to\Top$ a functor, then $\Map(\hocolim_{c\in\mathcal{C}}\mathcal{X}(c),Z)\simeq\holim_{c\in\mathcal{C}}\Map(\mathcal{X}(c),Z)$. In particular, it takes coproducts to products, and turns homotopy cofiber sequences into homotopy fiber sequences.
\end{prop}

\begin{proof}
We will indicate some of the ideas that go into showing a special case of this: if 
$$\xymatrix{
X_\emptyset \ar[r]\ar[d]& X_1\ar[d]\\
X_2\ar[r] & X_{12}
}
$$
is a homotopy pushout square, then for any space $Z$,
$$\xymatrix{
\Map(X_{12},Z)\ar[r]\ar[d]& \Map(X_{1},Z)\ar[d]\\
\Map(X_{2},Z)\ar[r] & \Map(X_{\emptyset},Z)
}
$$
is a homotopy pullback square. First, every homotopy pushout square admits an equivalence from a \emph{pushout square}; that is, a square of the form
$$\xymatrix{
X_\emptyset \ar[r]\ar[d]& X_1\ar[d]\\
X_2\ar[r] & X_{1}\cup_{X_\emptyset}X_{2}
}
$$
where $X_\emptyset\to X_i$ is a cofibration for each $i=1,2$. If we apply $\Map(-,Z)$ to this square, one checks by inspection that the resulting square is a pullback; that is, that $\Map(X_1\cup_{X_\emptyset}X_2,Z)=\Map(X_1,Z)\times_{\Map(X_\emptyset,Z)}\Map(X_2,Z)$. Hence the square
$$\xymatrix{
\Map(X_{12},Z)\ar[r]\ar[d]& \Map(X_{1},Z)\ar[d]\\
\Map(X_{2},Z)\ar[r] & \Map(X_{\emptyset},Z)
}
$$
is a pullback. The functor $\Map(-,Z)$ takes a cofibration $A\to X$ to a fibration $\Map(X,Z)\to\Map(A,Z)$, and so the pullback square is in fact a homotopy pullback.
\end{proof}

The proof of the general statement is in fact easier given an explicit description of the homotopy colimit of a functor $\mathcal{X}:\mathcal{C}\to\Top$ in terms of a under categories. See \cite{BK}.

\subsection{The Blakers-Massey Theorem}

We will only make statements for square diagrams, as those are the only types of cubical diagrams we have seriously studied in this work. All of what we say here has generalizations to higher dimensional cubes, and we refer the reader to \cite{CalcII} for details. We have already made use of the notion of a $k$-cartesian cube. Its dual notion, namely what it means for a cube to be $k$-cocartesian, is useful because the Blakers-Massey Theorem tells us how cartesian a $k$-cocartesian cube is.

\begin{defin}\label{D:totalcofiber}
For a $|T|$-cube $\mathcal{X}$, the \textsl{total homotopy cofiber} is the homotopy cofiber of the canonical map
$$
b(\mathcal{X}):\hocolim_{S\subsetneq T}\mathcal{X}(S).
$$
If $b(\mathcal{X})$ is $k$-connected, we say the cube is \textsl{$k$-cocartesian}, and if $k=\infty$, we say the cube is \textsl{homotopy cocartesian}.
\end{defin}

Thus a square 
$$\xymatrix{
X_\emptyset\ar[r]\ar[d]& X_1\ar[d]\\
X_2\ar[r] & X_{12}
}
$$
is homotopy cocartesian if the map $\hocolim(X_1\leftarrow X_\emptyset\rightarrow X_2)\to X_{12}$ is an equivalence. The homotopy colimit here is the double mapping cylinder.

\begin{thm}\label{T:BM}
Suppose the square
$$\xymatrix{
X_\emptyset\ar[r]\ar[d]& X_1\ar[d]\\
X_2\ar[r] & X_{12}
}
$$
is homotopy cocartesian, and that the maps $X_\emptyset\to X_i$ are $k_i$-connected. Then the square is $(k_1+k_2-1)$-cartesian.
\end{thm}

Here is a useful generalization.

\begin{thm}\label{T:BMgen}
Suppose the square
$$\xymatrix{
X_\emptyset\ar[r]\ar[d]& X_1\ar[d]\\
X_2\ar[r] & X_{12}
}
$$
is $k$-cocartesian, and that the maps $X_\emptyset\to X_i$ are $k_i$-connected. Then the square is $\min\{k-1,k_1+k_2-1\}$-cartesian.
\end{thm}

We will provide only a very bare sketch of \refT{BM}, if only to point out that one way of proving this uses disjunction techniques reminiscent of our arguments in \refS{disj}.

\begin{proof}
One can reduce to the case where $X_i$ is the union of $X_\emptyset$ with a $(k_i+1)$-cell $e^{k_i+1}$, and $X_{12}$ is the union of $X_\emptyset$ with both cells. We are therefore interested in the square
$$
\xymatrix{
X_\emptyset\ar[r]\ar[d]& X_\emptyset\cup e^{k_1+1}\ar[d]\\
X_\emptyset\cup e^{k_2+1}\ar[r] & X_\emptyset\cup e^{k_1+1}\cup e^{k_2+1}.
}
$$
For $i=1,2$, let $a_i$ be points in the interior of the two cells, and rewrite this square as 
$$
\xymatrix{
Y-\{a_1,a_2\}\ar[r]\ar[d]& Y-\{a_2\}\ar[d]\\
Y-\{a_1\}\ar[r] & Y.
}
$$
The claim is that the map
\begin{equation}\label{E:bmarg}
Y-\{a_1,a_2\}\to\holim(Y-\{a_2\}\rightarrow Y\leftarrow Y-\{a_1\})
\end{equation}
is $(k_1+k_2-1)$-connected. Recall the description of the codomain we gave just after \refD{holim}. A map $\gamma:S^k\to\holim(Y-\{a_2\}\rightarrow Y\leftarrow Y-\{a_1\})$ corresponds by adjointness to a map $\tilde{\gamma}:S^k\times I\to Y$ such that $\tilde{\gamma}(s,0)\neq a_2$ and $\tilde{\gamma}(s,1)\neq a_1$ for all $s\in S^k$. Extend $\tilde{\gamma}$ to $(-\epsilon, 1+\epsilon)$ to avoid talking about manifolds with corners. By a small homotopy, make $\tilde{\gamma}$ smooth near the $a_i$ and transverse to them (we may speak of smoothness because the interior of a cell has a smooth structure; transversality here means the $a_i$ are regular values of $\tilde{\gamma}$). Now consider the map $\Gamma: S^k\times \{(t_1,t_2)|t_i\neq t_2\}\to Y\times Y$ 
given by $(s,t_1,t_2)\mapsto (\tilde{\gamma}(t_1), \tilde{\gamma}(t_2))$. Again a small homotopy will make $a_1\times a_2$ a regular value of $\Gamma$. Note that $\Gamma^{-1}(a_1\times a_2) $ has dimension $k-k_1-k_2$ and hence will be empty if $k<k_1+k_2$. This means that for such $k$, $\tilde{\gamma}$ is homotopic to a map $\gamma'(s,t)$ which misses the $a_i$ for all $t$, and so $\gamma'$ lifts to a map $S^k\to Y-\{a_1,a_2\}$. Hence the map in \refE{bmarg} is surjective on homotopy groups for $k\leq k_1+k_2-1$. A similar argument establishes injectivity when $k<k_1+k_2-1$.
\end{proof}

\subsection{Acknowledgments}

The author thanks Jes\'us Gonz\'alez for helpful conversations during the preparation of this manuscript, and thanks Wellesley College for their hospitality.

\bibliographystyle{amsplain}

\bibliography{Bibliography}

\end{document}